\documentclass{amsart}

\usepackage{sty/preamble}

\title{Multifunctorial Inverse \texorpdfstring{$K$}{K}-Theory}

\authorinfoNJDY

\hypersetup{pdfauthor=\authors}

\date{15 July 2022}


\begin{document}

\begin{abstract}
  We show that Mandell's inverse $K$-theory functor is a categorically-enriched non-symmetric multifunctor.
In particular, it preserves algebraic structures parametrized by non-symmetric operads.
As applications, we describe how ring categories arise as the images of inverse $K$-theory.

\end{abstract}

\maketitle

\tableofcontents

\def\dprod{\displaystyle\prod}
\def\dsum{\displaystyle\sum}
\def\txsma{\textstyle\bigwedge}
\def\mtil{\til{m}}
\def\xtil{\til{x}}
\def\ytil{\til{y}}

\section{Introduction}\label{sec:introduction}

Segal's construction of $\Ga$-spaces from permutative categories \cite{segal} provides a functor
\[
\Kse \cn \permcat \to \Gasset,
\]
where $\permcat$ denotes the category of small permutative categories and $\Gasset$ denotes the category of $\Ga$-simplicial sets.
Stable equivalences of permutative categories are created by $\Kse$.
Work of Thomason \cite{thomason} shows that $\Kse$ induces an equivalence of stable homotopy categories.
For a detailed description of Segal $K$-theory, see \cite[Chapter~8]{cerberusIII}.

Mandell \cite[4.3]{mandell_inverseK} defines a homotopy inverse to $\Kse$ as a composite
\begin{equation}\label{eq:cScP}
  \Gasset \fto{\cS_*} \GaCat \fto{\cP} \permcatsus.
\end{equation}
We review the definition of $\cP$ in \cref{sec:mandellP} below and discuss $\cS_*$ after \cref{thm:invK_multifunctor}.
Here $\GaCat$ denotes the category of small $\Ga$-categories, and $\permcatsus$ denotes the subcategory of $\permcat$ consisting of strictly unital strong symmetric monoidal functors.
Following earlier work of Thomason \cite{thomason-model-str,cisinski-thomason-correction}, the functor $\cS_*$ induces a Quillen equivalence between the stable model structure on $\GasSet$ and a corresponding model structure on $\GaCat$.

We let $\Kinv = \cP \circ \cS_*$ denote the composite \cref{eq:cScP}, denoted
by $P$ in \cite[4.3]{mandell_inverseK}.
We use the term \emph{inverse $K$-theory} to refer to either $\Kinv$ or $\cP$, depending on the context.
The main results of \cite{mandell_inverseK} combine to show that $\Kinv$ induces an adjoint equivalence of stable homotopy categories.
This result provides a new proof and a sharpening of the stable equivalence due to Thomason \cite{thomason}.
\begin{theorem}[{\cite{mandell_inverseK}}]
  The functors $\Kse$ and $\Kinv$ induce an adjoint equivalence of stable homotopy categories.
\end{theorem}

In this paper we show that $\cP$ in \cref{eq:cScP} is a non-symmetric multifunctor in the $\Cat$-enriched sense. This result implies that $\cP$ preserves algebraic structures parametrized by non-symmetric operads such as the associative operad. The context for our main theorem is the following.
\begin{itemize}
\item The category $\GaCat$ is bicomplete and closed symmetric monoidal, with monoidal products given by Day convolution of pointed diagrams.  As a consequence, $\Gacat$ is a $\Cat$-enriched multicategory.
  We review these structures in \cref{sec:Ga-objs}.
\item There is a $\Cat$-enriched multicategory structure on $\permcatsus$, with operations given by multilinear functors and multilinear transformations.
  We review this in \cref{sec:permcatsu-multicat}.
\end{itemize} 
Now we state the main result.
\begin{theorem}\label{thm:invK_multifunctor}
The inverse $K$-theory functor
\[\GaCat \fto{\cP} \permcatsus\]
is a non-symmetric $\Cat$-enriched multifunctor.
\end{theorem}

Two remarks are in order.
\begin{enumerate}
\item The functor $\cS_*$ in \cref{eq:cScP} is \emph{not} a multifunctor even in the non-symmetric sense because it is the levelwise application of the composite functor
\[\begin{tikzcd}[column sep=large]
\sSet \ar{r}{\Sd^2} \ar[bend left=25]{rr}{\cS} & \sSet \ar{r}{h} & \Cat.
\end{tikzcd}\]
Here
\begin{itemize}
\item $\Sd$ is the subdivision functor, and
\item $h$ is the left adjoint of the nerve functor $N \cn \Cat \to \sSet$.
\end{itemize}
The subdivision functor $\Sd$ is not a monoidal functor, although it is oplax monoidal.
\item While \cref{thm:invK_multifunctor} shows that $\cP$ is a non-symmetric multifunctor in the $\Cat$-enriched sense, it is \emph{not} a multifunctor in the symmetric sense.  In \cref{remark:pseudo-symmetry} we explain that the symmetry diagram for $\cP$ commutes up to isomorphism, but does not strictly commute in general.
\end{enumerate}

As noted in \cite[pages~181--182]{elmendorf-mandell}, the Segal $\Gamma$-category functor
\[\begin{tikzcd}[column sep=large]
  \permcatsu \ar{r}{\Se{(-)}} & \Gacat
\end{tikzcd}\]
is \emph{not} a multifunctor, even ignoring the symmetry axiom, because its definition is incompatible with any potential pairing on a permutative category $\C$.
This limitation of $\Se{(-)}$, and hence also $\Kse$, is the main motivation for the work of Elmendorf-Mandell in \cite{elmendorf-mandell,elmendorf-mandell-perm}:
\begin{itemize}
\item replacing the indexing category $\Fskel$ in Segal $K$-theory with the more elaborate category $\Gskel$ and
\item developing the associated Elmendorf-Mandell $K$-theory
  \[\begin{tikzcd}[column sep=large]
  \permcatsu \ar{r}{\Kem} & \SymSp,
  \end{tikzcd}\]
  which is a simplicially-enriched multifunctor to symmetric spectra and is level equivalent to Segal $K$-theory.
\end{itemize} 
Therefore it is somewhat unexpected that the functor $\cP$, whose domain
$\Gacat$ uses the indexing category $\Fskel$, \emph{is} compatible with the multicategory structures of $\GaCat$ and $\permcatsus$. In this sense, Mandell's inverse $K$-theory is better
behaved than Segal $K$-theory.
As an application, we observe in \cref{sec:preservation} that $\cP$ sends associative monoids in $\Gacat$ to those in $\permcatsus$ (\Cref{Kinven}).

\subsection*{Outline}

\Cref{sec:general-background,sec:permcatsu-multicat,sec:Ga-objs} cover definitions and basic results that we will use.
Most of that material is well-known to experts, and we give references for more detailed treatment in the literature.
The definition of Mandell's inverse $K$-theory functor, $\cP$, is in  \Cref{sec:mandellP}.
The proof that $\cP$ is a non-symmetric multifunctor occupies the next three sections:
\begin{itemize}
\item The assignment on multimorphisms is given in \cref{sec:multimorphismassignment}.
\item The linearity constraints of each $\cP F$ are defined in \cref{sec:linearityconstraints}.
\item The non-symmetric enriched multifunctoriality axioms for $\cP$ are checked in \cref{sec:Pmultifunctor}.
\end{itemize}
In \cref{sec:preservation} we give an application with the associative operad and ring category structures.

\subsection*{Alternative Approach}

While this article was under review, the authors became aware of independent prior work of Elmendorf \cite{elmendorf-multi-inverse} that contains a multifunctoriality result similar to that of \cref{thm:invK_multifunctor}.
The approach of \cite{elmendorf-multi-inverse} is substantially different from our approach below, and is of independent interest.  Our approach here develops the $\Cat$-enrichment that is necessary for our applications in \cref{sec:preservation}.

\subsection*{Acknowledgment}

The authors would like to thank Tony Elmendorf and the referee for helpful comments on earlier versions of this paper.

\section{General Background}\label{sec:general-background}

In this section we outline the basic definitions, terminology, and notation to be used below.
We refer the reader to \cite{johnson-yau,cerberusIII} for a more detailed treatment.

\subsection*{Symmetric Monoidal Categories}

\begin{definition}
  A \emph{symmetric monoidal category} $(\C,\otimes,\tu,\al,\la,\rho,\xi)$ consists of
  \begin{itemize}
  \item a category $\C$;
  \item a functor $\otimes \cn \C \times \C \to \C$, which is called the \emph{monoidal product};
  \item an object $\tensorunit\in\C$, which is called the \emph{monoidal unit};
  \item a natural isomorphism
    \[\begin{tikzcd}[column sep=large]
    (X \otimes Y) \otimes Z \rar{\alpha_{X,Y,Z}}[swap]{\cong} &
    X\otimes (Y \otimes Z)
    \end{tikzcd}\]
    for all objects $X,Y,Z \in \C$, which is called the \emph{associativity isomorphism};
  \item natural isomorphisms
    \[
    \begin{tikzcd}
      \tensorunit \otimes X \rar{\lambda_X}[swap]{\cong} & X
    \end{tikzcd}
    \andspace
    \begin{tikzcd}
      X \otimes \tensorunit \rar{\rho_X}[swap]{\cong} & X
    \end{tikzcd}
    \]
    for all objects $X \in \C$, which are called the \emph{left unit isomorphism} and the \emph{right unit isomorphism}, respectively; and
  \item a natural isomorphism
    \[
    X \otimes Y \fto[\iso]{\xi_{X,Y}} Y \otimes X
    \]
    for objects $X,Y \in \C$, which is called the \emph{braiding} or \emph{symmetry isomorphism}.
  \end{itemize}

  The associativity and unit isomorphisms satisfy unity and pentagon
  axioms.  The braiding isomorphism satisfies further symmetry, unit,
  and hexagon axioms.  See \cite[1.1.1,1.1.23]{cerberusIII} for a complete
  description.
\end{definition}

\begin{definition}\label{def:monoidalfunctor}
  Suppose $(\C,\otimes,\tu^\C)$ and $(\D,\otimes,\tu^\D)$ are symmetric monoidal categories.  A \emph{symmetric monoidal functor} $(F,F^2,F^0)$ from $\C$ to $\D$, also called a
  \emph{symmetric lax monoidal functor}, consists of
  \begin{itemize}
  \item a functor
    \[
    F\cn \C \to \D,
    \]
  \item a monoidal constraint
    \[
    F^2_{X,Y} \cn FX \otimes FY \to F(X\otimes Y),
    \]
    and
  \item a unit constraint
    \[
    F^0 \cn \tu^\D \to F(\tu^\C)
    \]
  \end{itemize}
  satisfying axioms for associativity, unity, and symmetry.  See
  \cite[1.1.6,1.1.23]{cerberusIII} for a complete description.  We say that
  $(F,F^2,F^0)$ is
  \begin{itemize}
  \item \emph{unital} if $F^0$ is invertible;
  \item \emph{strictly unital} if $F^0$ is an identity;
  \item \emph{strong} if both $F^0$ and $F^2$ are invertible; and
  \item \emph{strict} if both $F^0$ and $F^2$ are identities.
  \end{itemize}

  A \emph{symmetric oplax monoidal functor} $(F,F^2,F^0)$ from $\C$ to
  $\D$ consists of
  \begin{itemize}
  \item a functor
    \[
    F\cn \C \to \D,
    \]
  \item an oplax monoidal constraint
    \[
    F^2_{X,Y} \cn F(X \otimes Y) \to FX \otimes FY,
    \]
    and
  \item an oplax unit constraint
    \[
    F^0 \cn F(\tu^\C) \to \tu^\D
    \]
  \end{itemize}
  satisfying oplax associativity, unity, and symmetry axioms.
\end{definition}

\begin{definition}\label{definition:permcat-u-su-sus}
  A \emph{permutative category} is a symmetric monoidal category whose
  associativity and unit morphisms are identities. We let $\PermCat$
  denote the 2-category of small permutative categories, symmetric
  monoidal functors, and monoidal natural transformations.  We also
  use the following locally-full sub 2-categories consisting of the
  same objects but restricting the 1-cells.
  \begin{itemize}
  \item $\PermCatu$ has 1-cells given by unital symmetric monoidal
    functors.
  \item $\PermCatsu$ has 1-cells given by strictly unital symmetric
    monoidal functors.
  \item $\permcatsus$ has 1-cells given by strictly unital strong
    symmetric monoidal functors.
  \item $\permcatst$ has 1-cells given by strict symmetric monoidal
    functors.
  \end{itemize}
  In each case the 2-cells are given by monoidal natural
  transformations.  We use the same notation for each of the
  underlying 1-categories.
\end{definition}

\subsection*{Multicategories}\label{sec:multicat}

We sketch the definitions of $\V$-enriched multicategories, multifunctors, and multinatural transformations.
Our two cases of interest will be $(\V,\otimes) = (\Set,\times)$ and $(\V,\otimes) =(\Cat,\times)$.
All of this material is covered in \cite[6.1]{cerberusIII} with complete details, additional context, and examples.
See \cite{bluemonster} and \cite{yau-hqft} for further development of the theory.

\begin{definition}\label{def:profile}
  Suppose $C$\label{notation:s-class} is a class.  
  \begin{enumerate}
  \item Denote by\label{notation:profs}
    \[\Prof(C) = \coprodover{n \geq 0}\ C^{\times n}\] 
    the class of finite ordered sequences of elements in $C$.
    An element in $\Prof(C)$ is called a \emph{$C$-profile}.  
  \item A typical $C$-profile of length $n=\len\angc$ is denoted by $\angc = (c_1, \ldots, c_n) \in C^{\times n}$\label{notation:us} or by $\ang{c_i}_i$ to indicate the indexing variable.
    The empty $C$-profile is denoted by $\ang{}$.
  \item We let $\oplus$ denote the concatenation of profiles, and note that $\oplus$ is an associative binary operation with unit given by the empty tuple $\ang{}$.
  \item An element in $\Prof(C)\times C$ is denoted as\label{notation:duc} $\IMMduc$ with $c'\in C$ and $\angc\in\Prof(C)$.
    \defmark
  \end{enumerate}
\end{definition}

\begin{definition}\label{def:enr-multicategory}
  Suppose $(\V,\otimes,\tu)$ is a symmetric monoidal category.  A
  \emph{$\V$-enriched multicategory} $(\M, \gamma, \operadunit)$
  consists of the following data.
  \begin{itemize}
  \item $\M$ is equipped with a class $\ObM$ of \emph{objects}.  We
    write $\Prof(\M)$ for $\Prof(\Ob\M)$.
  \item For $c'\in\ObM$ and $\angc=(c_1,\ldots,c_n)\in\ProfM$, $\M$ is
    equipped with an object of $\V$
    \[
    \M\IMMduc = \M\mmap{c'; c_1,\ldots,c_n} \in \V
    \]
    called the \emph{$n$-ary operation object} with \emph{input
    profile} $\angc$ and \emph{output} $c'$.
  \item For $\IMMduc \in \ProfMM$ as above and a permutation $\sigma
    \in \Sigma_n$, $\M$ is equipped with an isomorphism in $\V$
    \[\begin{tikzcd}
    \M\IMMduc \rar{\sigma}[swap]{\cong} & \M\IMMducsigma,
    \end{tikzcd}\]
    called the \emph{right action} or the \emph{symmetric group
    action}, in which
    \[
    \angc\sigma = (c_{\sigma(1)}, \ldots, c_{\sigma(n)})
    \]
    is the right permutation of $\angc$ by $\sigma$.
  \item For $c \in \ObM$, $\M$ is equipped with a
    morphism
    \[
    \operadunit_c\cn \tensorunit \to \M\IMMcc,
    \]
    called the \emph{$c$-colored unit}.
  \item For $c'' \in \ObM$, $\ang{c'} = (c'_1,\ldots,c'_n) \in
    \ProfM$, and $\ang{c_j} = (c_{j,1},\ldots,c_{j,k_j}) \in \ProfM$
    for each $j\in\{1,\ldots,n\}$, let $\angc = \oplus_j\ang{c_j} \in
    \ProfM$ be the concatenation of the $\ang{c_j}$.  Then $\M$ is
    equipped with a morphism in $\V$
  \begin{equation}\label{eq:enr-defn-gamma}
    \begin{tikzcd}
      \M\mmap{c'';\ang{c'}} \otimes
      \bigotimes\limits_{j=1}^n \M\mmap{c_j';\ang{c_j}}
      \rar{\gamma}
      &
      \M\mmap{c'';\ang{c}}
    \end{tikzcd}
  \end{equation}
  called the \emph{composition}. 
  \end{itemize}
  These data are required to satisfy axioms for the symmetric group
  action, associativity, unity, and two equivariance conditions.  See
  \cite[6.1.1]{cerberusIII} for a detailed description.
  A $\V$-enriched multicategory is
  \emph{small} if its class of objects is a set.
\end{definition}

\begin{definition}\label{definition:multicat-special}
  \
  \begin{enumerate}
  \item A $\V$-enriched multicategory with only one object is called a
    \emph{$\V$-enriched operad}. If $\M$ is a $\V$-enriched operad,
    then its object of $n$-ary operations is denoted by $\M_n \in \V$.
  \item A \emph{multicategory} is a $\Set$-enriched multicategory, where $(\Set,\times)$ is the symmetric monoidal category of sets and functions with the Cartesian product as the monoidal product.   An \emph{operad} is a $\Set$-enriched operad.
  \item Each small permutative category $(\C,\oplus)$ has an associated
    \emph{endomorphism multicategory} denoted $\End(\C)$.  The objects of
    $\End(\C)$ are those of $\C$ and the operations are given via the
    monoidal product in $\C$ as
    \[
    \End(\C)\mmap{Y;\ang{X}} = \C(\oplus_i X_i, Y).
    \]
An empty $\oplus$ means the monoidal unit in $\C$.
\item For an object $c$ in a $\V$-enriched multicategory $\M$, the \emph{$\V$-enriched endomorphism operad} $\End(c)$ has a single object $c$ and $n$-ary operation object
\[\End(c)_n = \M\scmap{\ang{c}_{i=1}^n;c} \forspace n \geq 0.\]
Its symmetric group action, unit, and composition are restrictions of those in $\M$.
  \item The \emph{terminal multicategory} $\Mterm$ has a single object
    and a single $n$-ary operation for each $n \ge 0$.
  \item The \emph{initial operad} $\Mtu$ has a single object and no
    non-identity operations.\dqed
  \end{enumerate}
\end{definition}

\begin{definition}\label{def:enr-multicategory-functor}
A \emph{$\V$-enriched multifunctor} $F : \M \to \N$ between $\V$-enriched multicategories $\M$ and $\N$ consists of the following data:
\begin{itemize}
\item an assignment \[F : \ObM \to \ObN,\] where $\ObM$ and $\ObN$ are the classes of objects of $\M$ and $\N$, respectively, and
\item for each $\mmap{c';\ang{c}} \in \ProfMM$ with $\angc=(c_1,\ldots,c_n)$, a morphism in $\V$
\[F : \M\mmap{c';\ang{c}} \to \N\mmap{Fc';F\ang{c}},\] where $F\angc=(Fc_1,\ldots,Fc_n)$.
\end{itemize}
These data satisfy axioms for preservation of the symmetric group
action, the colored units, and the composition.  See
\cite[6.1.10]{cerberusIII} for a detailed description.

A \emph{non-symmetric $\V$-enriched multifunctor} $F\cn \M \to \N$ is defined by the same data above---an object assignment and component $\V$-morphisms---but it is not required to preserve the symmetric group action of $\M$ and $\N$.  Thus a non-symmetric $\V$-enriched multifunctor is only required to preserve the colored units and composition.
\end{definition}

\begin{definition}\label{definition:enr-Palg}
  Suppose $\P$ is a $\V$-enriched operad and $\M$ is a $\V$-enriched
  multicategory.  A \emph{$\P$-algebra} in $\M$ is a pair
  \[
  (c,\theta)
  \]
  consisting of an object $c$ in $\M$ and a $\V$-enriched multifunctor
  \[
  \theta \cn \P \to \M
  \]
  that sends the single object of $\P$ to $c$.  Equivalently, $\theta$
  is a $\V$-enriched operad morphism
  \[
  \theta\cn \P \to \End(c).
  \]
  We say that $(c,\theta)$ is a \emph{non-symmetric} $\P$-algebra if
  \[
  \theta\cn \P \to \End(c)
  \]
  is a non-symmetric $\V$-enriched multifunctor.
\end{definition}

\begin{definition}\label{def:enr-multicat-natural-transformation}
Suppose $F,G : \M\to\N$ are $\V$-multifunctors as in \cref{def:enr-multicategory-functor}.  A \emph{$\V$-enriched multinatural transformation} $\alpha : F\to G$ consists of morphisms in $\V$
\[\alpha_c \cn \tu \to \N\mmap{Gc;Fc} \forspace c\in\ObM\] 
satisfying a $\V$-naturality axiom for each $\mmap{c';\ang{c}} \in
\ProfMM$.  See \cite[6.1.15]{cerberusIII} for a detailed description.
\end{definition}

\subsection*{Pointed Objects, Smash Products, and Pointed Homs}\label{sec:pointed}

\begin{definition}
  Suppose $(\C,\otimes,\tu,\Hom)$ is a complete and cocomplete symmetric monoidal
  closed category with terminal object $T$.  We let $\pC$ denote the category
  under $T$.  Its objects are morphisms
  \[
  \iota^X\cn T \to X \in \C,
  \]
  which are called \emph{pointed objects} with $\iota^X$ the
  \emph{basepoint} of $X$.  The morphisms of $\pC$ are called
  \emph{pointed morphisms} and are those morphisms of $\C$ that
  preserve the structure morphisms $\iota$.

  Moreover, we have the following constructions for pointed objects
  $X$ and $Y$.
  \begin{description}
  \item[Wedge] The \emph{wedge sum} $X \wed Y$ is the pushout in $\C$ of the span
  \[
  X \xleftarrow{\iota^X} T \fto{\iota^Y} Y.
  \]
  Equivalently, it is the coequalizer of the two induced morphisms
  \[
  \begin{tikzpicture}[x=30mm,y=20mm]
    \draw[0cell] 
    (0,0) node (ux) {T}
    (1,0) node (uy) {X \bincoprod Y}
    (1.8,0) node (w) {X \wed Y.}
    ;
    \draw[1cell] 
    (ux) edge[transform canvas={yshift=.8mm}] node {\iota^X} (uy)
    (ux) edge[transform canvas={yshift=-.8mm}] node['] {\iota^Y} (uy)
    (uy) edge[dashed] node {} (w)
    ;
  \end{tikzpicture}
  \]

  \item[Smash] The \emph{smash product} $X \sma Y$ is the following pushout in $\C$.
  \begin{equation}\label{eq:smash-pushout}
    \begin{tikzpicture}[x=50mm,y=15mm,vcenter]
      \draw[0cell=.9] 
      (0,0) node (a) {(X \otimes T) \bincoprod (T \otimes Y)}
      (1,0) node (b) {X \otimes Y}
      (0,-1) node (c) {T}
      (1,-1) node (d) {X \sma Y}
      ;
      \draw[1cell=.9] 
      (a) edge node {
        (1_X \otimes \iota_Y) \bincoprod (\iota_X \otimes 1_Y)
      } (b)
      (c) edge node {} (d)
      (a) edge node {} (c)
      (b) edge node {} (d)
      ;
    \end{tikzpicture}
  \end{equation}

  The \emph{smash unit} $E$ is defined by adjoining a disjoint
  basepoint to the unit of $\C$:
  \[
  E = \tu_+ = \tu \bincoprod T.
  \]

  \item[Pointed Hom] The \emph{pointed hom} $\pHom(X,Y)$ is the following pullback in $\C$.
  \begin{equation}\label{eq:pHom}
  \begin{tikzpicture}[x=30mm,y=15mm,vcenter]
    \draw[0cell] 
    (0,0) node (a) {\pHom(X,Y)}
    (1,0) node (b) {T}
    (0,-1) node (c) {\Hom(X,Y)}
    (1,-1) node (d) {\Hom(T,Y)}
    ;
    \draw[1cell] 
    (a) edge node {} (b)
    (c) edge node {} (d)
    (a) edge node {} (c)
    (b) edge node {} (d)
    ;
  \end{tikzpicture}
  \end{equation}
  The composite
  \[
  T \iso \Hom(X,T) \to \Hom(X,Y) \to \Hom(T,Y)
  \]
  induced by the structure morphisms for $X$ and $Y$ is equal to the vertical morphism in \cref{eq:pHom}, and therefore induces a canonical structure morphism $T \to \pHom(X,Y)$ making $\pHom(X,Y)$ a pointed object.\defmark
  \end{description}
\end{definition}

The coproduct of pointed objects is given by the wedge product.
Coequalizers of pointed objects are given by coequalizers in $\C$
equipped with the induced basepoint.  Similarly, products and
equalizers in $\pC$ are given in $\C$ and equipped with the basepoints
induced by these constructions.  Therefore $\pC$ is complete and
cocomplete if $\C$ is so.  The smash product and pointed hom provide a
symmetric monoidal closed structure with monoidal unit $E$.

\begin{theorem}[{\cite[4.20]{elmendorf-mandell-perm}, \cite[4.2.3]{cerberusIII}}]\label{theorem:pC-sm-closed}
  Suppose $(\C,\otimes,\tu,\Hom)$ is a complete and cocomplete
  symmetric monoidal closed category with terminal object $T$.  Then
  \[
  (\pC,\sma,E,\pHom)
  \]
  is a complete and cocomplete symmetric monoidal closed category.
\end{theorem}

We will use \cref{theorem:pC-sm-closed} in the following cases:
\begin{itemize}
\item $(\C,\otimes,\tu,T) = (\FinSet,\times,*,*)$, the category of
  finite sets, and
\item $(\C,\otimes,\tu,T) = (\Cat,\times,\boldone,\boldone)$, the
  category of small categories.
\end{itemize}

\begin{remark}\label{remark:multicat-smash}
  The smash product of small pointed multicategories is defined by taking $(\C,\otimes,\tu,T) = (\Multicat,\otimes,\Mtu,\Mterm)$, where $\otimes$ is the Boardman-Vogt tensor product.
  See \cite[5.6]{cerberusIII} for description of these products
  and their properties.
  Although we will not have use for it in this article, \cref{proposition:n-lin-equiv} gives a description of multilinear functors and transformations using the smash product of endomorphism multicategories.
\end{remark}

\section{The Multicategory of Small Permutative Categories}\label{sec:permcatsu-multicat}

In this section we give the definitions of multilinear functor between permutative categories and multilinear transformation from \cite{elmendorf-mandell} and \cite[6.5]{cerberusIII}.
These define the objects and morphisms for a $\Cat$-enriched multicategory also denoted $\permcatsu$ in \cref{definition:permcat-multicat} below.
Unless otherwise indicated, we generally use the following notation for permutative categories:
\begin{itemize}
\item $\oplus$ denotes the monoidal product.
\item $e$ denotes the monoidal unit.
\item $\xi$ denotes the braiding.
\end{itemize}

\begin{definition}\label{definition:tuple-subs}
  Suppose given permutative categories $\C_1,\ldots,\C_n$ and suppose $\ang{X} \in \txprod_j \C_j$ is a tuple of objects.
  \begin{itemize}
  \item For an object $X_i' \in \C_i$ with $i \in \{1,\ldots,n\}$, we
    let $\ang{X \compi X'_i}$ denote the tuple whose $j$th entry is
    $X_j$ for $j \ne i$, and whose $i$th entry is $X'_i$.
  \item Similarly, for objects $X'_i \in \C_i$ and $X_k' \in \C_k$
    with $i,k \in \{1,\ldots,n\}$ and $i \ne k$, we let $\ang{X \compi
      X_i' \compk X_k'}$ denote the tuple whose $j$th entry is that of
    $\ang{X \compi X_i'}$ for $j \ne k$ and whose $k$th entry is
    $X_k'$.\dqed
  \end{itemize}
\end{definition}

\begin{definition}[Multilinear Functors]\label{definition:multilinear-functors}
  Suppose $\C_1,\ldots,\C_n$, and $\D$ are permutative categories.  An
  \emph{$n$-linear functor} from $\ang{\C}$ to $\D$ is a functor
  \[
  F\cn \C_1 \times \cdots \times \C_n \to \D
  \]
  together with, for each $i\in
  \{1,\ldots,n\}$, a natural transformation $F^2_i$ called the
  \emph{ith linearity constraint} with components
  \[
  F^2_i\cn F\ang{X\compi X_i} \oplus F\ang{X\compi X_i'} \to
  F\ang{X\compi (X_i \oplus X_i')}
  \]
  for $\ang{X} \in \txprod_j \C_j$ and $X_i' \in \C_i$.  These data satisfy the following axioms.
  \begin{description}
  \item[Unity] If any $X_j = e$, the unit of $\C_j$, then $F\ang{X} = e$, the
    unit in $\D$.  Moreover, $F\ang{f \compj 1_e} = 1_e$ for any
    morphisms $f_i \in \C_i$ for $i \not= j$. 
  \item[Constraint Unity] If any $X_j = e$ or if $X_i'=e$, then $F^2_i$ is an identity
    morphism.
  \item[Constraint Associativity] The following diagram
    commutes for each $i\in \{1,\ldots,n\}$ and $\ang{X} \in \txprod_j
    \C_j$, with $X_i', X_i'' \in \C_i$.
    \begin{equation}\label{eq:ml-f2-assoc}
    \begin{tikzpicture}[x=55mm,y=15mm,vcenter]
      \draw[0cell=.85] 
      (0,0) node (a) {
        F\ang{X\compi X_i}
        \oplus F\ang{X\compi X_i'}
        \oplus F\ang{X\compi X_i''}
      }
      (1,0) node (b) {
        F\ang{X\compi X_i}
        \oplus F\ang{X\compi (X_i' \oplus X_i'')}
      }
      (0,-1) node (c) {
        F\ang{X\compi (X_i \oplus X_i')}
        \oplus F\ang{X\compi X_i''}
      }
      (1,-1) node (d) {
        F\ang{X\compi (X_i \oplus X_i' \oplus X_i'')}
      }
      ;
      \draw[1cell=.85] 
      (a) edge node {1 \oplus F^2_i} (b)
      (a) edge['] node {F^2_i \oplus 1} (c)
      (b) edge node {F^2_i} (d)
      (c) edge node {F^2_i} (d)
      ;
    \end{tikzpicture}
    \end{equation}
  \item[Constraint Symmetry] The following diagram commutes
    for each $i\in \{1,\ldots,n\}$ and $\ang{X} \in \txprod_j \C_j$, with
    $X_i' \in \C_i$.
    \begin{equation}\label{eq:ml-f2-symm}
    \begin{tikzpicture}[x=40mm,y=15mm,vcenter]
      \draw[0cell=.85] 
      (0,0) node (a) {
        F\ang{X\compi X_i}
        \oplus F\ang{X\compi X_i'}
      }
      (1,0) node (b) {
        F\ang{X\compi (X_i \oplus X_i')}
      }
      (0,-1) node (c) {
        F\ang{X\compi X_i'}
        \oplus F\ang{X\compi X_i}
      }
      (1,-1) node (d) {
        F\ang{X\compi (X_i' \oplus X_i)}
      }
      ;
      \draw[1cell=.85] 
      (a) edge node {F^2_i} (b)
      (a) edge['] node {\xi} (c)
      (b) edge node {F\ang{1 \compi \xi}} (d)
      (c) edge node {F^2_i} (d)
      ;
    \end{tikzpicture}
    \end{equation}
  \item[Constraint 2-By-2] The following diagram commutes
    for each
    \[
    i \not= k\in \{1,\ldots,n\}, 
    \quad
    \ang{X} \in \txprod_j \C_j,
    \quad
    X_i' \in \C_i,
    \andspace
    X_k' \in \C_k.
    \]
    \begin{equation}\label{eq:ml-f2-2by2}
    \begin{tikzpicture}[x=40mm,y=18mm,vcenter]
      \tikzset{0cell-nomm/.style={commutative diagrams/every diagram,
          every cell,nodes={scale=.8}
        }
      }
      \draw[0cell-nomm] 
      (0,0) node[align=left] (a) {$
        \phantom{\oplus}F\ang{X\compi X_i \compk X_k}
        \oplus F\ang{X\compi X_i' \compk X_k}$\\       
        $\oplus F\ang{X\compi X_i \compk X_k'}
        \oplus F\ang{X\compi X_i' \compk X_k'}
      $}
      (0,-1) node[align=left] (a') {$
        \phantom{\oplus}F\ang{X\compi X_i \compk X_k}
        \oplus F\ang{X\compi X_i \compk X_k'}$\\
        $\oplus F\ang{X\compi X_i' \compk X_k}
        \oplus F\ang{X\compi X_i' \compk X_k'}
      $}
      (.7,.75) node (b) {$
        F\ang{X\compi (X_i \oplus X_i') \compk X_k}
        \oplus F\ang{X\compi (X_i \oplus X_i') \compk X_k'} 
      $}
      (.7,-1.75) node (b') {$
        F\ang{X\compi X_i \compk (X_k \oplus X_k')}
        \oplus F\ang{X\compi X_i' \compk (X_k \oplus X_k')}
      $}
      (1.2,-.5) node (c) {$
        F\ang{X\compi (X_i \oplus X_i') \compk (X_k \oplus X_k')}
      $}
      ;
      \draw[1cell=.85] 
      (a) edge node[pos=.1] {F^2_i \oplus F^2_i} (b)
      (b) edge[transform canvas={shift={(.1,0)}}] node {F^2_k} (c)
      (a) edge['] node {1 \oplus \xi \oplus 1} (a')
      (a') edge['] node[pos=.1] {F^2_k \oplus F^2_k} (b')
      (b') edge[',transform canvas={shift={(.1,0)}}] node {F^2_i} (c)
      ;
    \end{tikzpicture}
    \end{equation}
  \end{description}

  A \emph{$0$-linear functor} is a choice of
  object in $\D$, regarded as a functor
  \[
  F \cn \boldone \to \D
  \]
  from the empty product.  We say that $F$ is a \emph{multilinear
    functor} if it is $n$-linear for some $n \ge 0$.  We say that a
  multilinear functor $F$ is
  \begin{itemize}
  \item \emph{strong} if $F^2_i$ is an isomorphism for each $i$ and
  \item \emph{strict} if $F^2_i$ is an identity for each $i$.\dqed
  \end{itemize}
\end{definition}

\begin{example}\label{explanation:1-linear-is-susm}
  In the context of \cref{definition:multilinear-functors} with $n=1$,
  comparing the axioms above with those of
  \cref{def:monoidalfunctor} shows that a 1-linear functor is
  precisely a strictly unital symmetric monoidal functor.
\end{example}

\begin{definition}[Multilinear Transformations]\label{explanation:multilinear-transformations}
  Suppose $\C_1,\ldots,\C_n$, and $\D$ are permutative categories.
  Suppose
  \[
  F,F'\cn \ang{\C} \to \D
  \]
  are $n$-linear functors.  An \emph{$n$-linear transformation}
  is a natural transformation of underlying functors
  \[
  \al \cn F \to F'
  \]
  that satisfies the following two \emph{multilinearity conditions}.
  \begin{enumerate}
  \item The diagram
  \begin{equation}\label{eq:monoidal-in-each-variable}
  \begin{tikzpicture}[x=45mm,y=15mm,vcenter]
      \draw[0cell=.85] 
      (0,0) node (a) {
        F\ang{X\compi X_i}
        \oplus F\ang{X\compi X_i'}
      }
      (1,0) node (b) {
        F\ang{X\compi (X_i \oplus X_i')}
      }
      (0,-1) node (c) {
        F'\ang{X\compi X_i}
        \oplus F'\ang{X\compi X_i'}
      }
      (1,-1) node (d) {
        F'\ang{X\compi (X_i \oplus X_i')}
      }
      ;
      \draw[1cell=.85] 
      (a) edge node {F^2_i} (b)
      (a) edge['] node {\al \oplus \al} (c)
      (b) edge node {\al} (d)
      (c) edge node {(F')^2_i} (d)
      ;
  \end{tikzpicture}
  \end{equation}
  commutes for each $i \in \{1,\ldots,n\}$ and $\ang{X} \in \txprod_j
  \C_j$ with $X_i' \in \C_i$.
  \item The component of $\al$ at a tuple $\ang{X}$ is an identity if
    any $X_i = e$.
  \end{enumerate}

  The two multilinearity conditions make $\al$ a monoidal natural transformation in each variable separately.
  We say that $\al$ is a \emph{multilinear transformation} if it is $n$-linear for some $n \ge 0$.
\end{definition}

\Cref{remark:multicat-smash} above reviews the smash product of
small pointed multicategories.  Although we will not have use for it in this
article, we point out that multilinearity can be defined equivalently
via multifunctors and multinatural transformations out of a smash
product of endomorphism multicategories.
\begin{proposition}[{\cite[6.5.10~and~6.5.13]{cerberusIII}}]\label{proposition:n-lin-equiv}
  For small permutative categories $\C_1,\ldots,\C_n$, and $\D$, the category of $n$-linear functors and $n$-linear transformations
  \[
  \ang{\C} \to \D
  \]
  is isomorphic to the category of multifunctors and multinatural transformations
  \begin{equation}\label{eq:multilin-via-sma}
    \txsma_{i=1}^n \End(\C_i) \to \End(\D).
  \end{equation}
\end{proposition}

\begin{definition}\label{definition:permcat-multicat}
  Let $\permcatsu$ denote the $\Cat$-enriched multicategory whose category of $n$-ary operations
  \[
  \permcatsu\mmap{\D;\ang{\C}}
  \]
  is the category of $n$-linear functors and $n$-linear transformations
  \[
  \ang{\C} \to \D.
  \]
  The $\Cat$-enriched multicategory axioms for $\permcatsu$ follow from \cref{proposition:n-lin-equiv} and the symmetric monoidal axioms of the smash product.
  Independently, a direct verification is given in \cite[6.6]{cerberusIII}.

  We let $\permcatsus$ denote the $\Cat$-enriched sub-multicategory
  whose $n$-ary operations consist of \emph{strong} $n$-linear functors and
  $n$-linear transformations.
\end{definition}

\begin{remark}
  The multicategory structure on $\permcatsu$ is not the endomorphism multicategory of a symmetric monoidal structure.
  For example, the unit for the smash product of multicategories is not a permutative category.
  See \cite[5.7.23 and 10.2.17]{cerberusIII} for further discussion of this point.
\end{remark}

\section{Symmetric Monoidal Closed Categories of \texorpdfstring{$\Ga$}{Gamma}-Objects}\label{sec:Ga-objs}

In this section we recall the definitions of $\Ga$-categories and $\Ga$-simplicial sets.

\begin{definition}\label{definition:FinSet-Fskel}
  Let $\pFinSet$ denote the category of pointed finite sets and pointed functions.
  Let $\Fskel$ denote the full subcategory of $\pFinSet$ whose objects are $\ord{n} = \{0,\ldots,n\}$ with basepoint $0$, for natural numbers $n \ge 0$.
  The pointed finite set $\ord{0}$ is both initial and terminal in $\Fskel$.
\end{definition}

\begin{definition}\label{definition:lexord}
  Suppose $\ord{m}$ and $\ord{n}$ are objects of $\Fskel$.
  The \emph{lexicographic order} for $\ord{m}\sma\ord{n}$ is the
  bijection of pointed finite sets
  \[
  L\cn\ord{m} \sma \ord{n} \iso \ord{mn}
  \]
  given by
  \[
  L(x,y) = 
  \begin{cases}
    0 & \text{if } x=0 \orspace y=0,\\
    n(x-1) + y & \text{if } x>0 \andspace y>0.
  \end{cases}
  \]
\end{definition}

Using the lexicographic order, the smash product of pointed finite
sets induces a monoidal product on $\Fskel$ that we also denote
$\sma$.  Elementary algebra with the formula for $L$ shows that this
product is strictly associative and unital with strict monoidal unit
$\ord{1}$.  We state this as the following result.
\begin{proposition}\label{prop:Fskel-perm}
  With the lexicographic ordering of smash products,
  $(\Fskel,\sma,\ord{1})$ is a permutative category.
\end{proposition}

\begin{definition}\label{def:pointedcategories}\
\begin{enumerate}
\item A \emph{pointed category} is a pair $(\C,*)$ consisting of a category $\C$ and a distinguished object $*$, called the \emph{basepoint}.
\item For pointed categories $(\C,*)$ and $(\D,*)$, a \emph{pointed functor}
\[\begin{tikzcd}[column sep=large]
(\C,*) \ar{r}{F} & (\D,*)
\end{tikzcd}\]
is a functor $F \cn \C \to \D$ that preserves the basepoint: $F(*) = *$.  Identity and composite pointed functors are defined by the underlying functors 
\item Suppose $F,G \cn (\C,*) \to (\D,*)$ are pointed functors.  A \emph{pointed natural transformation} $\theta \cn F \to G$ is a natural transformation with identity basepoint component:
\[\theta_* = 1_* \cn F(*) = * \to * = G(*) \in \D.\]
Identity, horizontal composite, and vertical composite pointed natural transformations are defined by the underlying natural transformations.\defmark
\end{enumerate}
\end{definition}

\begin{definition}\label{definition:gamma-obj}
  Suppose $(\C,*)$ is a pointed category with $*$ terminal in $\C$.  A
  \emph{$\Ga$-object} in $\C$ is a pointed functor
  \[
  X\cn (\Fskel,\ord{0}) \to (\C,*). 
  \]
  The category of $\Ga$-objects in $\C$, denoted $\Ga\mh\C$, is
  \[
  \pCat\big((\Fskel,\ord{0}),(\C,*)\big),
  \]
  the category of pointed functors from $(\Fskel,\ord{0})$ to $(\C,*)$
  and pointed natural transformations.
\end{definition}

If $\C$ is a symmetric monoidal category, then functors from another
symmetric monoidal category, such as $\Fskel$, have a symmetric
monoidal product given by Day convolution.  To show that this
preserves basepoints, and therefore induces a symmetric monoidal
product on $\Ga\mh\C$, we use enrichment over the category of pointed
sets, $\pSet$, via the following result.
\begin{lemma}[{\cite[4.3.5]{cerberusIII}}]\label{lemma:zCat-enr}
  Suppose $(\B,*)$ and $(\C,*)$ are pointed categories in which
  the basepoints are both initial and terminal.
  \begin{enumerate}
  \item\label{it:zCat-enr} Taking zero morphisms as basepoints of their hom sets, both
    $\B$ and $\C$ are enriched categories over $\pSet$.
  \item\label{it:zCat-equiv} There is an equivalence of categories
    \begin{equation}\label{eq:zCat-equiv}
      \pCat((\B,*), (\C,*)) \hty (\pSet\Cat)(\B,\C),
    \end{equation}
    where the right hand side denotes the category of $\pSet$-enriched
    functors and natural transformations.
  \end{enumerate}
\end{lemma}

Our applications below take $(\B,*)$ to be $(\Fskel,\ord{0})$ and take $(\C,*)$
to be either $(\pCat,\boldone)$ or $(\psSet,*)$.

\begin{definition}\label{definition:Fhat}
  Suppose $(\C,\otimes,\tu,*)$ is a symmetric monoidal closed category
  that is complete and cocomplete with chosen terminal object $*$.
  \begin{enumerate}
  \item For $\ord{a}$ and $\ord{b}$ in $\Fskel$, we use the
    notation
    \begin{equation}\label{eq:Gpunc}
      \Fpunc(\ord{a}\,,\ord{b}) = \big(\Fskel(\ord{a}\,,\ord{b})\big)^\punc
    \end{equation}
    for the subset of nonzero morphisms.
  \item We let $\Fhat$ denote
    $\Fskel$ equipped with the \emph{pointed unitary enrichment} over $\pC$
    defined by
    \begin{equation}\label{eq:GhatE}
    \Fhat\big(\ord{a}\,, \ord{b}\big) =
    \bigvee_{\Fpunc(\ord{a}\scs\ord{b})} S^0
    \end{equation}
    where $S^0 = \tu \bincoprod *$ is the monoidal unit of $(\pC,\sma)$.
    The empty wedge is the terminal object $*$ of $\C$.
  \item The monoidal product of $\pC$-categories is denoted $\sma$
    and, in particular, $\Fhat \sma \Fhat$ has objects given by pairs
    and hom objects given by
    \[
    \big(\Fhat\sma\Fhat\big)\big((\ord{a}\,,\ord{a'})\scs(\ord{b}\,,\ord{b'})\big)
    = \Fhat(\ord{a}\,,\ord{b})\sma \Fhat(\ord{a'}\,,\ord{b'})
    \]
    for natural numbers $a$, $a'$, $b$, and $b'$.\dqed
  \end{enumerate}
\end{definition}

\begin{definition}\label{definition:convolution-map-Fskel}
  Suppose $(\C,\otimes,\tu,[,],*)$ is a symmetric monoidal closed category
  that is complete and cocomplete with chosen terminal object $*$.
  Suppose given $\Ga$-objects $X$ and $Y$.  We
  define a mapping object $\pMap(X,Y)$ in $\pC$ and $\Ga$-objects
  \[
  X \sma Y, \quad \pHom(X,Y), \andspace \dtu
  \]
  as follows.

  \begin{itemize}
  \item The \emph{Day convolution product} of $X$ and $Y$ is
  \begin{equation}\label{eq:FotimesXY}
    X \sma Y =
    \ecint^{(\ord{a}\scs\ord{b}) \in \Fhat\sma\Fhat}
    \bigvee_{\Fpunc(\ord{a} \sma \ord{b}\scs -)}
    \big(X\ord{a} \sma Y \ord{b}\big).
  \end{equation}

  \item The \emph{hom diagram} for $X$ and $Y$ is
  \begin{align}\label{eq:FHomXY}
    \pHom(X,Y) & = \ecint_{(\ord{b}\scs\ord{c}) \in \Fhat\sma\Fhat}
    \;\prod_{\Fpunc(- \sma \ord{b}\scs \ord{c})}
    \big[X\ord{b}\scs Y\ord{c}]_*\\
    & \iso \ecint_{\ord{b} \in \Fhat}\, [X{\ord{b}}\scs Y({- \sma \ord{b}})]_*
  \end{align}
  where $[-,-]_*$ denotes the pointed internal hom of $(\pC,\sma)$.

  \item The \emph{mapping object} for $X$ and $Y$ is
  \begin{equation}\label{eq:FMapXY}
    \pMap(X,Y) = \ecint_{\ord{b} \in \Fhat}\, [X{\ord{b}}\scs Y{\ord{b}}]_* =
    \big(\pHom(X,Y)\big)\ord{1},
  \end{equation}
  where the second equality holds because $\ord{1}$ is a strict unit in
  $\Fskel$.
  
  \item The \emph{monoidal unit diagram} is
  \begin{equation}\label{eq:Funit}
    \dtu = \Fhat\big(\ord{1}\,,-\big) = \bigvee_{\Fpunc(\ord{1}\scs-)} S^0.
  \end{equation}

\item The \emph{braiding}
  \begin{equation}\label{eq:Fxi}
    X \sma Y \to Y \sma X
  \end{equation}
  is the natural transformation induced by the symmetry of the smash product
  \[
  X\ord{a} \sma Y\ord{b} \fto{\xi} Y\ord{b} \sma X\ord{a},
  \]
  the symmetry in $\Fskel$, and the universal properties of coproducts and coends.
  \end{itemize}
  These define a symmetric monoidal closed structure for $\Ga\mh\C$ by
  \cite[4.3.37]{cerberusIII}.

  Moreover, evaluation at $\ord{1}$ defines a symmetric monoidal functor
  \begin{equation}\label{eq:pC-GstarC-adj}
    \begin{tikzpicture}[x=30mm,y=25mm,vcenter]
      \draw[0cell] 
      (0,0) node (x) {L_{\ord{1}}\cn \pC}
      (1,0) node (y) {\Ga\mh\C \cn \ev_{\ord{1}}.}
      ;
      \draw[1cell] 
      (x) edge[transform canvas={yshift=1.4mm}] node {} (y) 
      (y) edge[transform canvas={yshift=0mm}] node {} (x) 
      ;
    \end{tikzpicture}
  \end{equation}
  The left adjoint $L_{\ord{1}}$ is defined for $A \in \pC$ by composition
  \[
  A \sma \dtu \cn \Fhat \fto{\dtu} \pC \fto{A \sma -} \pC.\dqed
  \]
\end{definition}

\begin{theorem}[{\cite[4.3.37~and~9.2.18]{cerberusIII}}]\label{corollary:Ga-C-smenr}
  Suppose $(\C,\otimes,\tu,[,],*)$ is a symmetric monoidal closed category that is complete and cocomplete with chosen terminal object $*$.
  Then $\Ga\mh\C$ is a complete and cocomplete symmetric monoidal closed category that is enriched, tensored, and cotensored over $\pC$.
\end{theorem}


\section{Mandell's Inverse \texorpdfstring{$K$}{K}-Theory Functor \texorpdfstring{$\cP$}{P}}
\label{sec:mandellP}

In this section we review Mandell's inverse $K$-theory functor $\cP$ \cite{mandell_inverseK}.
For an integer $n \geq 0$, 
\begin{itemize}
\item $\ufs{n}$ denotes the unpointed finite set $\{1,\ldots,n\}$ with $\ufs{0} = \emptyset$, and
\item $\ord{n}$ denotes the pointed finite set $\{0,1,\ldots,n\}$ with basepoint 0.
\end{itemize}
For a pointed finite set $a$ with basepoint $*$, we write 
\[a^\punc = a \setminus \{*\}\]
for the unpointed finite set obtained from $a$ by removing the basepoint.  Note that $\ord{n}^\punc = \ufs{n}$.

First recall from \cite[Section 4]{mandell_inverseK} that, for a $\Gamma$-category $X$, the small permutative category $\cP X$ is defined as the Grothendieck construction
\begin{equation}\label{mandellPX}
\cP X = \int_{\cA} AX.
\end{equation}

\subsection*{The Category $\cA$}

The indexing category $\cA$ has, as objects, finite sequences of positive integers.  A morphism 
\begin{equation}\label{Amorphismphi}
\begin{tikzcd}[column sep=large]
m = (m_1, \ldots, m_p) \ar{r}{\phi} & (n_1, \ldots, n_q) = n \in \cA
\end{tikzcd}
\end{equation}
is a map of unpointed finite sets
\[\begin{tikzcd}[column sep=large]
\coprod\limits_{i=1}^p \ufs{m_i} \ar{r}{\phi} & \coprod\limits_{j=1}^q \ufs{n_j}
\end{tikzcd}\]
such that the preimage of each $\ufs{n_j}$ is either empty or contained in a single $\ufs{m_i}$, with the index $i$ determined by the index $j$.  Note that, however, each $\ufs{m_i}$ may hit several different $\ufs{n_j}$.

\subsection*{The Functor $AX$}
On objects the functor 
\[\begin{tikzcd}[column sep=large]
\cA \ar{r}{AX} & \Cat
\end{tikzcd}\]
is the assignment, for $m_1,\ldots,m_p > 0$,  
\[(AX)(m_1, \ldots, m_p) = 
\begin{cases}
\prod\limits_{i=1}^p X\ord{m_i} & \text{if $p>0$ and}\\
X\ord{0} = \boldone & \text{if $p=0$.}
\end{cases}\]
For a morphism $\phi \cn m \to n$ in $\cA$ as in \cref{Amorphismphi}, the functor
\begin{equation}\label{phistarAphi}
\begin{tikzcd}[column sep=2.5cm]
(AX)(m) \ar{r}{\phi_* \,=\, (AX)(\phi)} & (AX)(n)
\end{tikzcd}
\end{equation}
is the identity functor of the terminal category $\boldone$ if $q = 0$, which forces $p=0$.  If $p=0$ and $q>0$, then the composite of $\phi_*$ with the $j$-th coordinate projection, for $j \in \{1,\ldots,q\}$, is the following composite.
\[\begin{tikzcd}[column sep=large]
\boldone \ar[equal]{d} \ar{r}{\phi_*} & (AX)(n) = \prod\limits_{j=1}^q X\ord{n_j} \ar[shorten <=-1ex]{d}{\text{project}}\\
X\ord{0} \ar{r} & X\ord{n_j}
\end{tikzcd}\]
The bottom horizontal arrow is induced by the pointed map $\ord{0} \to \ord{n_j}$ in $\Fskel$.  Its image is the basepoint of $X\ord{n_j}$. 

For $p,q > 0$ and $j \in \{1,\ldots,q\}$, if $\phi^{\inv}\big(\ufs{n_j}\big) = \emptyset$, then the composite of $\phi_*$ with the $j$-th coordinate projection is the following composite, whose image is the basepoint of $X\ord{n_j}$.
\[\begin{tikzcd}[column sep=large]
(AX)(m) = \prod\limits_{i=1}^p X\ord{m_i} \ar[shorten <=-1ex]{d} \ar{r}{\phi_*} 
& (AX)(n) = \prod\limits_{j=1}^q X\ord{n_j} \ar[shorten <=-1ex]{d}{\text{project}}\\
\boldone = X\ord{0} \ar{r} & X\ord{n_j}
\end{tikzcd}\]
If $\phi^{\inv}\big(\ufs{n_j}\big) \not= \emptyset$, then there is a unique index $i \in \{1,\ldots,p\}$ such that 
\[\emptyset \not= \phi^{\inv}\big(\ufs{n_j}\big) \subset \ufs{m_i}.\]
Define the map of pointed finite sets
\begin{equation}\label{phiijpointed}
\begin{tikzcd}[column sep=large]
\ord{m_i} \ar{r}{\phi_{i,j}} & \ord{n_j} \in \Fskel
\end{tikzcd}
\end{equation}
by
\[\phi_{i,j}(x) = \begin{cases}
\phi(x) & \text{if $x \in \phi^{\inv}\big(\ufs{n_j}\big) \subset \ufs{m_i}$ and}\\
0 & \text{if $x \in \ord{m_i} \setminus \phi^{\inv}\big(\ufs{n_j}\big)$}.
\end{cases}\]
The composite of $\phi_*$ with the $j$-th coordinate projection is the following composite.
\begin{equation}\label{xofphiij}
\begin{tikzcd}[column sep=large]
(AX)(m) = \prod\limits_{i=1}^p X\ord{m_i} \ar[shorten <=-1ex]{d}[swap]{\text{project}} \ar{r}{\phi_*} 
& (AX)(n) = \prod\limits_{j=1}^q X\ord{n_j} \ar[shorten <=-1ex]{d}{\text{project}}\\
X\ord{m_i} \ar{r}{(\phi_{i,j})_* \,=\, X\phi_{i,j}} & X\ord{n_j}
\end{tikzcd}
\end{equation}
\begin{remark}
  In \cite[page 777, line 2]{mandell_inverseK} the description of the map $\phi_{i,j}$ is slightly incorrect.
  The correct definition of $\phi_{i,j}$, given above, is from \cite[Notation~5.8]{gjo1}.
\end{remark}

\subsection*{The Category $\cP X$}
An object in the Grothendieck construction $\cP X = \int_{\cA} AX$ is a pair $(m,x)$ with
\begin{itemize}
\item $m = (m_1,\ldots,m_p)$ an object in $\cA$ and
\item $x = (x_1, \ldots, x_p)$ an object in $(AX)(m)$ with each $x_i \in X\ord{m_i}$.
\end{itemize}
A morphism 
\begin{equation}\label{pofxmorphism}
  \begin{tikzcd}[column sep=large]
    (m,x) \ar{r}{(\phi, f)} & (n,y) \in \cP X
  \end{tikzcd}
\end{equation}
consists of
\begin{itemize}
\item a morphism $\phi \cn m \to n$ in $\cA$ as in \cref{Amorphismphi} and 
\item a morphism $f \cn \phi_*(x) \to y \in (AX)(n)$ with $\phi_* = (AX)\phi$ as in \cref{phistarAphi}.
\end{itemize}
Composition is defined as
\[(\psi,g) \circ (\phi,f) = \big(\psi \phi, g \circ (\psi_* f)\big).\]
The identity morphism of an object $(m,x)$ is the pair $(1_m,1_x)$ of identity morphisms.  

For a $\Ga$-category morphism $F \cn X \to Y$, the functor 
\[\cP F \cn \cP X \to \cP Y\]
is defined on objects by
\[(\cP F)(m,x) = (m,Fx)\]
with $Fx = (Fx_1, \ldots, Fx_p)$ and similarly on morphisms.  

Moreover, $\cP X$ is a permutative category with monoidal product on objects given by concatenation in each variable, as in
\begin{equation}\label{PXmonoidalprod}
(m,x) \Box (n,y) = \big((m,n), (x,y)\big) \in \cP X
\end{equation}
with
\[\begin{split}
(m,n) &= \big(m_1, \ldots, m_p, n_1, \ldots, n_q\big) \in \cA\\
(x,y) &= \big(x_1, \ldots, x_p, y_1, \ldots, y_q\big) \in (AX)(m,n).
\end{split}\]
The monoidal product on morphisms is defined similarly, using the fact that $\cA$ is a permutative category under concatenation of sequences with $()$ as the monoidal unit.  The pair 
\begin{equation}\label{PXmonoidalunit}
\big((),*\big) \in \cP X
\end{equation}
consisting of
\begin{itemize}
\item the empty tuple $()$ in $\cA$ and
\item the unique object $* \in X\ord{0} = \boldone$
\end{itemize} 
is defined as the strict unit for $\Box$.  The braiding
\[(m,x) \Box (n,y) \fto{\xi = (\tau, 1)} (n,y) \Box (m,x) \in \cP X\]
consists of
\begin{itemize}
\item the isomorphism
\[\tau \cn (m,n) \fto{\iso} (n,m) \in \cA\]
given by the block swapping
\[\begin{tikzcd}
\left(\coprod_{i=1}^p \ufs{m_i}\right) \coprod \left(\coprod_{j=1}^q \ufs{n_j}\right) \ar{r}{\tau} & \left(\coprod_{j=1}^q \ufs{n_j}\right) \coprod \left(\coprod_{i=1}^p \ufs{m_i}\right)
\end{tikzcd}\]
of unpointed finite sets and
\item the identity morphism of $(y,x)$.
\end{itemize}

As stated in \cite[4.5]{mandell_inverseK}, the above definitions define a functor
\[\begin{tikzcd}[column sep=large]
\Gacat \ar{r}{\cP} & \permcatst
\end{tikzcd}\]
with $\permcatst$ the category of small permutative categories and strict symmetric monoidal functors.  We also let $\cP$ denote the composite with the inclusion
\[
\permcatst \hookrightarrow \permcatsus
\]
of strict symmetric monoidal functors among strictly unital strong symmetric monoidal functors.

\section{Multimorphism Assignment}
\label{sec:multimorphismassignment}

To show that
\[\begin{tikzcd}[column sep=large]
\Gacat \ar{r}{\cP} & \permcatsus
\end{tikzcd}\]
is a non-symmetric $\Cat$-enriched multifunctor, suppose $X_1, \ldots, X_k$, and $Z$ are $\Gamma$-categories for $k \geq 0$, and
\[(X_1, \ldots, X_k) \to Z\]
is a $k$-morphism in the multicategory $\Gacat$.  This is a $\Gamma$-category morphism
\begin{equation}\label{smashxfz}
\begin{tikzcd}[column sep=large]
\bigwedge\limits_{i=1}^k X_i \ar{r}{F} & Z
\end{tikzcd}
\end{equation}
with $\sma$ being the Day convolution of pointed diagrams from \cref{definition:convolution-map-Fskel}.
Unless otherwise specified, an iterated monoidal product, such as $\bigwedge_{i=1}^k$, is  assumed to be left normalized.
We need to construct its image $k$-morphism
\[
(\cP X_1, \ldots, \cP X_k) \to \cP Z
\]
in the multicategory $\permcatsus$.
This is a strong $k$-linear functor
\begin{equation}\label{prodpxpfpz}
  \begin{tikzcd}[column sep=2.8cm]
    \prod\limits_{i=1}^k \cP X_i \ar{r}{\big(\cP F,\, \{(\cP F)^2_i\}_{i=1}^k \big)} & \cP Z
  \end{tikzcd}
\end{equation}
between small permutative categories.
There are two cases depending on whether $k=0$ or $k > 0$.

\subsection*{The case $k=0$}

If $k=0$, then the given $\Gamma$-category morphism in \cref{smashxfz} is
\begin{equation}\label{JFZ}
\begin{tikzcd}[column sep=large]
\ftu = \bigvee\limits_{\Fpunc(\ord{1}, -)} S^0 \ar{r}{F} & Z
\end{tikzcd}
\end{equation}
with
\begin{itemize}
\item $\ftu \in \Gacat$ the monoidal unit with respect to $\sma$,
\item $\Fpunc(\ord{1}, \ord{n})$ the set of nonzero morphisms $\ord{1} \to \ord{n}$ in $\Fskel$ for each $\ord{n}$, and
\item $S^0 = \{0,1\}$ the two-object discrete category with basepoint 0.
\end{itemize}
For each $n \geq 0$, there is a bijection 
\[\Fpunc(\ord{1}, \ord{n}) \iso \ord{n}^\punc = \{1,\ldots,n\}.\]
There is a pointed functor
\[\begin{tikzcd}[column sep=large]
\ftu\ord{n} = \bigvee\limits_{\Fpunc(\ord{1}, \ord{n})} S^0 \iso \ord{n} \ar{r}{F_n} & Z\ord{n}
\end{tikzcd}\]
that is natural in $\ord{n} \in \Fskel$.  For each $j \in\{1,\ldots,n\}$, the pointed map
\[\begin{tikzcd}[column sep=large,row sep=0ex,
/tikz/column 1/.append style={anchor=base east},
/tikz/column 2/.append style={anchor=base west}]
\ord{1} \ar{r} & \ord{n} \in \Fskel\\
1 \ar[mapsto]{r} & j
\end{tikzcd}\]
yields a commutative diagram of pointed functors
\[\begin{tikzcd}[column sep=large]
\ord{1} \iso \ftu\ord{1} \ar{d} \ar{r}{F_1} & Z\ord{1} \ar{d}\\
\ord{n} \iso \ftu\ord{n} \ar{r}{F_n} & Z\ord{n}
\end{tikzcd}\]
that sends $1 \in \ftu\ord{1}$ to $F_1(1) \in Z\ord{1}$ and then to $F_n(j) \in Z\ord{n}$.  So $F$ is completely determined by $F_1(1)$.
 
The desired 0-linear functor to $\cP Z$ in \cref{prodpxpfpz} is a choice of an object in $\cP Z$.  We define the image of the 0-morphism $F$ in \cref{smashxfz} as the object
\begin{equation}\label{foneonePZ}
\big((1) \in \cA, F_1(1) \in (AZ)(1) = Z\ord{1}\big) \in \cP Z.
\end{equation}
This finishes the definition of the 0-linear functor $\cP F$ in \cref{prodpxpfpz} when $k=0$.

\subsection*{Object assignment for $k>0$}

For $k > 0$, the given $\Gamma$-category morphism $F$ in \cref{smashxfz} consists of component pointed functors
\begin{equation}\label{smashxfzr}
\begin{tikzcd}[column sep=large]
\Big(\bigwedge\limits_{i=1}^k X_i\Big)\ord{r} 
= \dint^{(\ord{p_1},\ldots,\ord{p_k}) \in \Fhat^{\sma k}} \mkern-20mu \bigvee\limits_{\Fpunc\left(\sma_i \ord{p_i},\, \ord{r}\right)} \bigwedge\limits_{i=1}^k \big(X_i\ord{p_i}\big) \ar{r}{F_r} & Z\ord{r}
\end{tikzcd}
\end{equation}
for $\ord{r} \in \Fskel$ that are natural in $\ord{r}$.  In the $\pCat$-coend above, an empty wedge is the terminal category $\boldone$.  We will sometimes omit the subscript $r$ in $F_r$.  

Equivalently, by the universal property of the $\pCat$-coend in \cref{smashxfzr}, $F$ is a \emph{pointed} natural transformation between pointed functors (\cref{def:pointedcategories}) as follows.
\[\begin{tikzpicture}
\def\u{-1.3}
\draw[0cell]
(0,0) node (a) {\Fskel^{\sma k}}
(a)+(3.5,0) node (b) {\pCat^{\sma k}}
(a)+(0,\u) node (c) {\Fskel}
(b)+(0,\u) node (d) {\pCat}
;
\draw[1cell]
(a) edge node {X_1 \overline{\sma} \cdots \overline{\sma} X_k} (b)
(b) edge node {\sma} (d)
(a) edge node[swap] {\sma} (c)
(c) edge node[swap] {Z} (d)
;
\draw[2cell]
node[between=a and d at .42, shift={(0,0)}, rotate=-130, 2label={above,F}] {\Rightarrow}
;
\end{tikzpicture}\]
In the top horizontal arrow, we use $\overline{\sma}$ to denote the smash product of pointed functors, which is different from the Day convolution in \cref{eq:FotimesXY}.  Unpacking the above pointed natural transformation, $F$ is uniquely determined by component pointed functors
\begin{equation}\label{Fcompnent}
\begin{tikzcd}[column sep=large]
\bigwedge_{i=1}^k X_i \ord{p_i} \ar{r}{F} & Z\left(\ord{p_1 \cdots p_k}\right)
\end{tikzcd}
\end{equation}
for $(\ord{p_1},\ldots,\ord{p_k}) \in \Fskel^{\sma k}$ that are compatible with morphisms in $\Fskel^{\sma k}$.  The composite pointed functor \cref{Fxonekjonekdef} below is of this form.

To define the $k$-linear functor $\cP F$ in \cref{prodpxpfpz}, first we define its assignment on objects.  Suppose given, for each $i \in \{1,\ldots,k\}$, an object
\begin{equation}\label{mixi}
(m^i,x^i) \in \cP X_i
\end{equation}
with
\[\begin{split}
m^i &= (m^i_1,\ldots,m^i_{r_i}) \in \cA\\
x^i &= (x^i_1,\ldots,x^i_{r_i}) \in (AX_i)(m^i) = \prod_{j=1}^{r_i} X_i\ord{m^i_j}.  
\end{split}\]
If $r_i = 0$ for some $i$, then
\[(m^i,x^i) = \big((),*\big) \in \cP X_i\]
is the monoidal unit, and we define the image object
\begin{equation}\label{pfmonoidalunit}
  (\cP F)\big((m^1,x^1), \ldots , (m^k,x^k)\big) = \big((),*\big) \in \cP Z,
\end{equation}
the monoidal unit in $\cP Z$, as part of the definition of $\cP F$ in \cref{prodpxpfpz}.
The definition \cref{pfmonoidalunit} is forced by the unity axiom of a multilinear functor in \cref{definition:multilinear-functors}.

Suppose $r_i > 0$ for each $i \in \{1,\ldots,k\}$.  Given $j_i \in \{1,\ldots,r_i\}$ for each $i$, we first define the positive integer
\begin{equation}\label{monekjonek}
m^{1,\ldots,k}_{j_1,\ldots,j_k} = \prod_{i=1}^k m^i_{j_i}.
\end{equation}
Letting each $j_i$ run through $\{1,\ldots,r_i\}$, we define the object
\begin{equation}\label{monek}
m^{1,\ldots,k} = \Big\{ \cdots \, \Big\{ m^{1,\ldots,k}_{j_1,\ldots,j_k} \Big\}_{j_1=1}^{r_1} \,\cdots\, \Big\}_{j_k=1}^{r_k} \in \cA
\end{equation}
with length $r_1\cdots r_k$.  In other words, the displayed integer $m^{1,\ldots,k}_{j_1,\ldots,j_k}$ is entry number
\[j_1 + \sum_{i=2}^k \big[ r_1 \cdots \, r_{i-1} (j_i - 1) \big] \]
in the sequence $m^{1,\ldots,k}$.  To understand the definitions \cref{monekjonek,monek}, we arrange the objects $m^1, \ldots, m^k \in \cA$ as follows.
\begin{equation}\label{eq:monek-rows}\begin{split}
m^1 &= \big( m^1_1, \ldots, m^1_{r_1} \big)\\
\vdots &\\
m^k &= \big( m^k_1, \ldots, m^k_{r_k} \big)
\end{split}\end{equation}
The positive integer $m^{1,\ldots,k}_{j_1,\ldots,j_k}$ is a vertical product in this arrangement, taking the $j_i$-th entry in the $i$-th row for each $i \in \{1,\ldots,k\}$.  In the object $m^{1,\ldots,k} \in \cA$, these positive integers are ordered reverse lexicographically by first comparing $j_k$, then $j_{k-1}$, and so forth, with $j_1$ compared last.

Similarly, define the objects
\begin{equation}\label{xonekjonek}
x^{1,\ldots,k}_{j_1,\ldots,j_k} = \big(x^1_{j_1}, \ldots, x^k_{j_k}\big) \in \bigwedge_{i=1}^k \Big(X_i \ord{m^i_{j_i}}\Big).
\end{equation}
Define the object
\begin{equation}\label{Fxonekjonek}
F\big(x^{1,\ldots,k}_{j_1,\ldots,j_k}\big) \in Z\big(\ord{m^{1,\ldots,k}_{j_1,\ldots,j_k}}\big)
\end{equation}
as the image of $x^{1,\ldots,k}_{j_1,\ldots,j_k}$ under the following composite of pointed functors.
\begin{equation}\label{Fxonekjonekdef}
\begin{tikzcd}[column sep=large,cells={nodes={scale=.85}}]
\bigwedge\limits_{i=1}^k \Big(X_i \ord{m^i_{j_i}}\Big) \ar{d}[swap]{\omega} &\\
\Big(\bigwedge\limits_{i=1}^k X_i\Big)\big(\ord{m^{1,\ldots,k}_{j_1,\ldots,j_k}}\big) 
= \dint^{(\ord{p_1},\ldots,\ord{p_k}) \in \Fhat^{\sma k}} \mkern-30mu \bigvee\limits_{\Fpunc\big(\sma_i \ord{p_i},\, \ord{m^{1,\ldots,k}_{j_1,\ldots,j_k}}\big)} \bigwedge\limits_{i=1}^k \left(X_i\ord{p_i}\right) \ar{r}{F} & Z\big(\ord{m^{1,\ldots,k}_{j_1,\ldots,j_k}}\big)
\end{tikzcd}
\end{equation}
In \cref{Fxonekjonekdef}, $\omega$ is the composite of (i) the coend structure morphism for the index
\[(\ord{p_1},\ldots,\ord{p_k}) = \big(\ord{m^1_{j_1}},\ldots,\ord{m^k_{j_k}}\big) \in \Fhat^{\sma k}\] 
and (ii) the wedge factor inclusion corresponding to the identity morphism of
\begin{equation}\label{smashmiji}
\bigwedge_{i=1}^k \ord{m^i_{j_i}} = \ord{m^{1,\ldots,k}_{j_1,\ldots,j_k}} \in \Fskel
\end{equation}
in the wedge index.  Alternatively, the composite in \cref{Fxonekjonekdef} is a component of the form \cref{Fcompnent}.  With entries from \cref{Fxonekjonek}, define the object
\begin{equation}\label{Fxonek}
\begin{split}
F(x^{1,\ldots,k}) &= \Big\{ \cdots \, \Big\{ F\big(x^{1,\ldots,k}_{j_1,\ldots,j_k}\big) \Big\}_{j_1=1}^{r_1} \,\cdots\, \Big\}_{j_k=1}^{r_k}\\
& \in \prod_{j_k=1}^{r_k} \cdots \prod_{j_1=1}^{r_1} Z\big(\ord{m^{1,\ldots,k}_{j_1,\ldots,j_k}}\big) = (AZ)(m^{1,\ldots,k}).
\end{split}
\end{equation} 
Using \cref{monek,Fxonek}, we define the image object
\begin{equation}\label{pfmx}
(\cP F)\big((m^1,x^1), \ldots , (m^k,x^k)\big) = \big(m^{1,\ldots,k}, F(x^{1,\ldots,k})\big) \in \cP Z.
\end{equation}
We can actually regard the image object \cref{pfmonoidalunit}, when $r_i = 0$ for some $i$, as a special case of \cref{pfmx} as follows: 
\begin{itemize}
\item $m^{1,\ldots,k} \in \cA$ has length $r_1\cdots r_k = 0$, so it is the empty sequence $()$.  
\item $F(x^{1,\ldots,k}) \in (AZ)() = Z\ord{0} = \boldone$ is the unique object.
\end{itemize} 
This finishes the definition of the object assignment of $\cP F$ in \cref{prodpxpfpz} when $k>0$.  In \cref{rk:redefinition} we explain in detail that we can also use the lexicographic ordering consistently instead of the reverse lexicographic ordering in the construction of $\cP F$.

\subsection*{Morphism assignment for $k>0$}

To define $\cP F$ in \cref{prodpxpfpz} on morphisms, suppose given, for each $i \in \{1,\ldots,k\}$, a morphism
\begin{equation}\label{phiifi}
\begin{tikzcd}[column sep=large]
(m^i,x^i) \ar{r}{(\phi^i,f^i)} & (n^i,y^i) \in \cP X_i
\end{tikzcd}
\end{equation}
with $(m^i,x^i)$ as in \cref{mixi}, objects
\begin{equation}\label{niyiobjects}
\begin{split}
n^i &= (n^i_1, \ldots , n^i_{s_i}) \in \cA\\
y^i &= (y^i_1, \ldots , y^i_{s_i}) \in (AX_i)(n^i) = \prod_{\ell=1}^{s_i} X_i\ord{n^i_\ell},
\end{split}
\end{equation}
and morphisms
\begin{equation}\label{phiifimorphisms}
\begin{tikzcd}[column sep=large,row sep=0ex,
/tikz/column 1/.append style={anchor=base east},
/tikz/column 2/.append style={anchor=base west}]
m^i \ar{r}{\phi^i} & n^i \in \cA\\
\phi^i_*(x^i) \ar{r}{f^i} & y^i \in (AX_i)(n^i).
\end{tikzcd}
\end{equation}
We want to construct a corresponding morphism 
\[\Big(\phi \cn m^{1,\ldots,k} \to n^{1,\ldots,k}, f \cn \phi_* F(x^{1,\ldots,k}) \to F(y^{1,\ldots,k}) \Big) \in \cP Z\]
as in \cref{pofxmorphism}.  

First we define the morphism in $\cA$
\begin{equation}\label{phimoneknonek}
\begin{tikzcd}[column sep=small,row sep=small]
m^{1,\ldots,k} \ar[equal]{d} \ar{r}{\phi} & n^{1,\ldots,k} \ar[equal]{d} \\
\Big\{ \cdots \, \Big\{ m^{1,\ldots,k}_{j_1,\ldots,j_k} \Big\}_{j_1=1}^{r_1} \,\cdots\, \Big\}_{j_k=1}^{r_k} & \Big\{ \cdots \, \Big\{ n^{1,\ldots,k}_{\ell_1,\ldots,\ell_k} \Big\}_{\ell_1=1}^{s_1} \,\cdots\, \Big\}_{\ell_k=1}^{s_k}
\end{tikzcd}
\end{equation}
with $m^{1,\ldots,k}_{j_1,\ldots,j_k}$ as in \cref{monekjonek} and similarly
\[n^{1,\ldots,k}_{\ell_1,\ldots,\ell_k} = \prod_{i=1}^k n^i_{\ell_i}
\forspace 1 \leq \ell_i \leq s_i,\] 
as the following map of unpointed finite sets.
\begin{equation}\label{phionephik}
\begin{tikzcd}[column sep=2.3cm]
\coprod\limits_{j_k=1}^{r_k} \cdots \coprod\limits_{j_1=1}^{r_1} \ufs{m^{1,\ldots,k}_{j_1,\ldots,j_k}} \ar{r}{\phi \,=\, (\phi^i)_{i=1}^k} & 
\coprod\limits_{\ell_k=1}^{s_k} \cdots \coprod\limits_{\ell_1=1}^{s_1} \ufs{n^{1,\ldots,k}_{\ell_1,\ldots,\ell_k}} 
\end{tikzcd}
\end{equation}
More precisely, for each $k$-tuple of elements
\[a = (a_1,\ldots,a_k) \in \prod_{i=1}^k \ufs{m^i_{j_i}} \iso \ufs{m^{1,\ldots,k}_{j_1,\ldots,j_k}}\]
with each $a_i \in \ufs{m^i_{j_i}}$, the given map of unpointed finite sets
\[\begin{tikzcd}[column sep=large]
\coprod\limits_{j_i=1}^{r_i} \ufs{m^i_{j_i}} \ar{r}{\phi^i} & \coprod\limits_{\ell_i=1}^{s_i} \ufs{n^i_{\ell_i}} 
\end{tikzcd}\]
sends $a_i$ to an element
\begin{equation}\label{phiiai}
\phi^i(a_i) \in \ufs{n^i_{\ell_i}}
\end{equation}
for some unique index $\ell_i \in \{1,\ldots,s_i\}$.  This yields the $k$-tuple of elements
\[\phi(a) = \big(\phi^1(a_1), \ldots , \phi^k(a_k) \big) \in \prod_{i=1}^k \ufs{n^i_{\ell_i}} \iso \ufs{n^{1,\ldots,k}_{\ell_1,\ldots,\ell_k}},\]
which is in one of the coproduct summands in the codomain of $\phi$ in \cref{phionephik}.  Since each $\phi^i \cn m^i \to n^i$ is a morphism in $\cA$, the condition 
\[\emptyset \not= (\phi^i)^{\inv}\big( \ufs{n^i_{\ell_i}} \big) \bigsubset \ufs{m^i_{j_i}} \forspace 1 \leq i \leq k\]
implies the condition
\begin{equation}\label{nonekelloneknotempty}
\emptyset \not= \phi^{\inv}\Big( \ufs{n^{1,\ldots,k}_{\ell_1,\ldots,\ell_k}} \Big) 
\bigsubset \ufs{m^{1,\ldots,k}_{j_1,\ldots,j_k}}.
\end{equation}
This ensures that $\phi$ in \cref{phimoneknonek} is a morphism in $\cA$.

The other component of the desired morphism in $\cP Z$ is a morphism
\begin{equation}\label{fphistarfxtofy}
\begin{tikzcd}[column sep=large]
\phi_* F(x^{1,\ldots,k}) \ar{r}{f} & F(y^{1,\ldots,k}) \in (AZ)(n^{1,\ldots,k}).
\end{tikzcd}
\end{equation}
To define $f$, first note that in the setting of \cref{phiiai}, there is an induced map of pointed finite sets
\begin{equation}\label{phijielli}
\begin{tikzcd}[column sep=large]
\ord{m^i_{j_i}} \ar{r}{\phi^i_{j_i, \ell_i}} & \ord{n^i_{\ell_i}} \in \Fskel
\end{tikzcd}
\end{equation}
as in \cref{phiijpointed}.  The morphism $\phi \cn m^{1,\ldots,k} \to n^{1,\ldots,k}$ in $\cA$ in \cref{phimoneknonek} induces a functor
\begin{equation}\label{AZphi}
\begin{tikzcd}[column sep=large, row sep=small]
(AZ)(m^{1,\ldots,k}) \ar[equal]{d} \ar{r}{\phi_* \,=\, (AZ)\phi} & (AZ)(n^{1,\ldots,k}) \ar[equal]{d}\\
\prod\limits_{j_k=1}^{r_k} \cdots \prod\limits_{j_1=1}^{r_1} Z \Big( \ord{m^{1,\ldots,k}_{j_1,\ldots,j_k}} \Big) & 
\prod\limits_{\ell_k=1}^{s_k} \cdots \prod\limits_{\ell_1=1}^{s_1} Z \Big( \ord{n^{1,\ldots,k}_{\ell_1,\ldots,\ell_k}} \Big) 
\end{tikzcd}
\end{equation}
as in \cref{phistarAphi}.  In the nontrivial cases, the coordinates of $\phi_*$ are induced by the morphisms $\phi^i_{j_i,\ell_i}$ in \cref{phijielli}.  More precisely, for a coordinate that satisfies the condition \cref{nonekelloneknotempty}, $\phi_*$ in \cref{AZphi} is given by the following composite.
\begin{equation}\label{AZphicoordinates}
\begin{tikzcd}[column sep=large]
\prod\limits_{j_k=1}^{r_k} \cdots \prod\limits_{j_1=1}^{r_1} Z \Big( \ord{m^{1,\ldots,k}_{j_1,\ldots,j_k}} \Big) \ar[shorten <=-1ex]{d}[swap]{\text{project}} \ar{r}{\phi_*} 
& \prod\limits_{\ell_k=1}^{s_k} \cdots \prod\limits_{\ell_1=1}^{s_1} Z \Big( \ord{n^{1,\ldots,k}_{\ell_1,\ldots,\ell_k}} \Big)  \ar[shorten <=-1ex]{d}{\text{project}}\\
Z \Big( \ord{m^{1,\ldots,k}_{j_1,\ldots,j_k}} \Big) \ar{r}{Z\big( \bigwedge_{i=1}^k \phi^i_{j_i,\ell_i} \big)} 
& Z \Big( \ord{n^{1,\ldots,k}_{\ell_1,\ldots,\ell_k}} \Big)
\end{tikzcd}
\end{equation}
The bottom horizontal functor is the image under $Z$ of the morphism
\[\begin{tikzcd}[column sep=2.5cm]
\ord{m^{1,\ldots,k}_{j_1,\ldots,j_k}} = \bigwedge\limits_{i=1}^k \ord{m^i_{j_i}} \ar{r}{\bigwedge_{i=1}^k \phi^i_{j_i,\ell_i}} 
& \bigwedge\limits_{i=1}^k \ord{n^i_{\ell_i}} = \ord{n^{1,\ldots,k}_{\ell_1,\ldots,\ell_k}} \in \Fskel,
\end{tikzcd}\]
which uses the definition \cref{smashmiji} of the smash product in $\Fskel$.

The given morphism $f^i \cn \phi^i_*(x^i) \to y^i$ in $(AX_i)(n^i)$ in \cref{phiifimorphisms} consists of component morphisms
\begin{equation}\label{fiellitoyielli}
\begin{tikzcd}[column sep=large]
\left(\phi^i_{j_i,\ell_i}\right)_* \big(x^i_{j_i}\big) \ar{r}{f^i_{\ell_i}} & y^i_{\ell_i} \in X_i\ord{n^i_{\ell_i}}
\end{tikzcd} \ifspace \emptyset \not= (\phi^i)^{\inv}\big(\ufs{n^i_{\ell_i}}\big) \bigsubset \ufs{m^i_{j_i}},
\end{equation}
with the functor
\[\left(\phi^i_{j_i,\ell_i}\right)_* = X_i\big(\phi^i_{j_i,\ell_i}\big) \cn X_i\ord{m^i_{j_i}} \to X_i\ord{n^i_{\ell_i}}.\]
The $\ell_i$-th entry of $f^i$ is a morphism 
\begin{equation}\label{startoyielli}
\begin{tikzcd}[column sep=large]
* \ar{r}{f^i_{\ell_i}} & y^i_{\ell_i} \in X_i\ord{n^i_{\ell_i}}
\end{tikzcd}
\ifspace (\phi^i)^{\inv}\big(\ufs{n^i_{\ell_i}}\big) = \emptyset.
\end{equation}

For a factor of $(AZ)(n^{1,\ldots,k})$ that satisfies the condition \cref{nonekelloneknotempty}, the composite pointed functor
\begin{equation}\label{smashxinitozn}
\begin{tikzcd}[column sep=huge]
\bigwedge\limits_{i=1}^k \Big( X_i \ord{n^i_{\ell_i}} \Big) \ar{r}{F \circ \omega} 
& Z \Big( \ord{n^{1,\ldots,k}_{\ell_1,\ldots,\ell_k}} \Big),
\end{tikzcd}
\end{equation}
defined as in \cref{Fxonekjonekdef}, sends the morphism
\begin{equation}\label{fielliinsmashx}
\big\{ f^i_{\ell_i} \big\}_{i=1}^k \in \bigwedge\limits_{i=1}^k \Big( X_i \ord{n^i_{\ell_i}} \Big),
\end{equation}
with each component of the form \cref{fiellitoyielli}, to a morphism
\begin{equation}\label{phistarfxtofy}
\begin{tikzcd}[column sep=huge]
F \Big\{ \big(\phi^i_{j_i,\ell_i}\big)_* \big(x^i_{j_i}\big) \Big\}_{i=1}^k \ar{r}{f^{1,\ldots,k}_{\ell_1,\ldots,\ell_k}} 
& F\big( y^{1,\ldots,k}_{\ell_1,\ldots,\ell_k} \big) \in Z \Big( \ord{n^{1,\ldots,k}_{\ell_1,\ldots,\ell_k}} \Big).
\end{tikzcd}
\end{equation}
Since $F$ in \cref{smashxfzr} is natural in the morphisms in $\Fskel$, the domain in \cref{phistarfxtofy} is also given by
\[F \Big\{ \big(\phi^i_{j_i,\ell_i}\big)_* \big(x^i_{j_i}\big) \Big\}_{i=1}^k 
= \Big( Z\big( \txsma_{i=1}^k \phi^i_{j_i,\ell_i} \big) \Big) \Big( F \big( x^{1,\ldots,k}_{j_1,\ldots,j_k} \big) \Big).\]
In this case, the $Z \big( \ord{n^{1,\ldots,k}_{\ell_1,\ldots,\ell_k}} \big)$ component of the desired morphism $f$ in \cref{fphistarfxtofy} is defined as the morphism in \cref{phistarfxtofy}.

For a factor of $(AZ)(n^{1,\ldots,k})$ that does \emph{not} satisfy \cref{nonekelloneknotempty}, we have
\[(\phi^i)^{\inv} \big( \ufs{n^i_{\ell_i}} \big) = \emptyset\]
for at least one index $i \in \{1,\ldots,k\}$.  The pointed functor $F \circ \omega$ in \cref{smashxinitozn} sends the morphism $\big\{ f^i_{\ell_i} \big\}_{i=1}^k$ in \cref{fielliinsmashx}, with at least one component of the form \cref{startoyielli}, to a morphism
\begin{equation}\label{fstartofy}
\begin{tikzcd}[column sep=huge]
* \ar{r}{f^{1,\ldots,k}_{\ell_1,\ldots,\ell_k}} 
& F\big( y^{1,\ldots,k}_{\ell_1,\ldots,\ell_k} \big) \in Z \Big( \ord{n^{1,\ldots,k}_{\ell_1,\ldots,\ell_k}} \Big).
\end{tikzcd}
\end{equation}
In \cref{fstartofy} the domain is the basepoint $*$ because
\begin{itemize}
\item $F \circ \omega$ is a pointed functor and
\item at least one entry in the input object is the basepoint, namely, $* \in X_i\ord{n^i_{\ell_i}}$, which yields the basepoint in $\bigwedge_{i=1}^k \big( X_i \ord{n^i_{\ell_i}} \big)$.
\end{itemize} 
In this case, the $Z \big( \ord{n^{1,\ldots,k}_{\ell_1,\ldots,\ell_k}} \big)$ component of the desired morphism $f$ in \cref{fphistarfxtofy} is defined as the morphism in \cref{fstartofy}.

In summary, given the morphisms \cref{phiifi} in $\cP X_i$ for $i \in \{1,\ldots,k\}$, we define the morphism
\begin{equation}\label{PFmorphism}
(\cP F)\big( (\phi^1,f^1), \ldots, (\phi^k,f^k) \big) = (\phi, f) \in \cP Z
\end{equation}
with
\begin{itemize}
\item $\phi \cn m^{1,\ldots,k} \to n^{1,\ldots,k} \in \cA$ as in \cref{phimoneknonek} and
\item $f \cn \phi_* F(x^{1,\ldots,k}) \to F(y^{1,\ldots,k})$ in \cref{fphistarfxtofy} having components \cref{phistarfxtofy,fstartofy}.
\end{itemize} 
This finishes the definition of the morphism assignment of $\cP F$ in \cref{prodpxpfpz} for $k > 0$.  The functoriality of $\cP F$ follows from the definitions \cref{phionephik,phistarfxtofy,fstartofy}, and the functoriality of $F \circ \omega$ in \cref{smashxinitozn}.

\section{Linearity Constraints}
\label{sec:linearityconstraints}

To make $\cP F$ in \cref{prodpxpfpz} for $k > 0$ into a strong $k$-linear functor, next we construct its $k$ linearity constraints $\{(\cP F)^2_i\}_{i=1}^k$ as in \cref{definition:multilinear-functors}.
Suppose given
\begin{itemize}
\item objects $(m^i,x^i) \in \cP X_i$ for $i \in \{1,\ldots,k\}$ as in \cref{mixi}, 
\item an index $b \in \{1,\ldots,k\}$, and
\item an object 
\[(\mhat^b, \xhat^b) \in \cP X_b\]
with
\[\begin{split}
\mhat^b &= \big( \mhat^b_1, \ldots, \mhat^b_{\rhat_b} \big) \in \cA\\
\xhat^b &= \big( \xhat^b_1, \ldots, \xhat^b_{\rhat_b} \big) \in (AX_b)(\mhat^b) = \prod_{j=1}^{\rhat_b} X_b\ord{\mhat^b_j}.
\end{split}\]
\end{itemize} 
With the notation
\begin{equation}\label{angmx}
\begin{split}
\ang{m,x} &= \big((m^1,x^1), \ldots, (m^k,x^k)\big) \in \prod_{i=1}^k \cP X_i\\
\ang{m,x} \circ_b (\mhat^b, \xhat^b) 
&= \big((m^1,x^1), \ldots, \underbrace{(\mhat^b, \xhat^b)}_{\text{$b$-th entry}}, \ldots, (m^k,x^k)\big), 
\end{split}
\end{equation}
the corresponding $b$-th linearity constraint of $\cP F$,
\begin{equation}\label{PFtwob}
\begin{tikzcd}[column sep=large]
(\cP F)\ang{m,x} \Box (\cP F)\big( \ang{m,x} \circ_b (\mhat^b, \xhat^b) \big) \ar{d}{(\cP F)^2_b}\\ 
(\cP F)\Big( \ang{m,x} \circ_b \big( (m^b,x^b) \Box (\mhat^b, \xhat^b) \big) \Big), 
\end{tikzcd}
\end{equation}
is the isomorphism in $\cP Z$ defined as follows.

\subsection*{The Domain}

In the domain of $(\cP F)^2_b$ in \cref{PFtwob}, the factor 
\[(\cP F)\ang{m,x} = \big( m^{1,\ldots,k}, F(x^{1,\ldots,k}) \big) \in \cP Z\]
is defined in \cref{pfmx}.  Replacing $(m^b,x^b)$ with $(\mhat^b,\xhat^b)$ yields the other factor,
\begin{equation}\label{PFmxcompb}
(\cP F)\big( \ang{m,x} \circ_b (\mhat^b, \xhat^b) \big) 
= \big( \mhat^{1,\ldots,k}, F(\xhat^{1,\ldots,k}) \big) \in \cP Z,
\end{equation}
with the following notation adapted from \cref{monekjonek,monek,xonekjonek,Fxonekjonek,Fxonek}.
\begin{equation}\label{mhatonekjonek}
\begin{split}
\mhat^{1,\ldots,k}_{j_1,\ldots,j_k} 
&= m^1_{j_1} \cdots \, \mhat^b_{j_b} \cdots \, m^k_{j_k} \qquad (\text{with $1 \leq j_b \leq \rhat_b$})\\
\mhat^{1,\ldots,k} 
&= \Big\{ \cdots \, \Big\{ \cdots \, \Big\{ \mhat^{1,\ldots,k}_{j_1,\ldots,j_k} \Big\}_{j_1=1}^{r_1} \,\cdots\, \Big\}_{j_b=1}^{\rhat_b} \,\cdots\, \Big\}_{j_k=1}^{r_k} \in \cA\\
\xhat^{1,\ldots,k}_{j_1,\ldots,j_k} 
&= \big(x^1_{j_1}, \ldots, \xhat^b_{j_b}, \ldots, x^k_{j_k}\big)\\
&\in \big(X_1 \ord{m^1_{j_1}}\big) \sma \cdots \sma \big(X_b \ord{\mhat^b_{j_b}}\big) \sma \cdots \sma \big(X_k \ord{m^k_{j_k}}\big)\\
F(\xhat^{1,\ldots,k}) 
&= \Big\{ \cdots \, \Big\{ \cdots \, \Big\{ F\big(\xhat^{1,\ldots,k}_{j_1,\ldots,j_k}\big) \Big\}_{j_1=1}^{r_1} \,\cdots\, \Big\}_{j_b=1}^{\rhat_b} \,\cdots\, \Big\}_{j_k=1}^{r_k}\\
& \in \prod_{j_k=1}^{r_k} \cdots \prod_{j_b=1}^{\rhat_b} \cdots \prod_{j_1=1}^{r_1} Z\big(\ord{\mhat^{1,\ldots,k}_{j_1,\ldots,j_k}}\big) = (AZ)(\mhat^{1,\ldots,k})
\end{split}
\end{equation}
In each of the definitions of $\mhat^{1,\ldots,k}_{j_1,\ldots,j_k}$ and $\xhat^{1,\ldots,k}_{j_1,\ldots,j_k}$ above, the entry involving $\widehat{?}^b_{j_b}$ appears in the $b$-th spot.  The monoidal product $\Box$ in the domain of $(\cP F)^2_b$ in \cref{PFtwob} is taken in $\cP Z$ as in \cref{PXmonoidalprod}.  Thus the domain of $(\cP F)^2_b$ in \cref{PFtwob} is the object
\begin{equation}\label{PFtwobdomain}
\begin{split}
&\big( m^{1,\ldots,k}, F(x^{1,\ldots,k}) \big) \Box \big( \mhat^{1,\ldots,k}, F(\xhat^{1,\ldots,k}) \big)\\
&= \Big( \big( m^{1,\ldots,k}, \mhat^{1,\ldots,k} \big), \big( F(x^{1,\ldots,k}), F(\xhat^{1,\ldots,k})\big) \Big) \in \cP Z
\end{split}
\end{equation}
with second entry the object $\big( F(x^{1,\ldots,k}), F(\xhat^{1,\ldots,k})\big)$ in
\[\scalebox{.95}{$
(AZ) \big( m^{1,\ldots,k}, \mhat^{1,\ldots,k} \big)
= \bigg[ \dprod_{j_k=1}^{r_k} \cdots \dprod_{j_1=1}^{r_1} Z\big(\ord{m^{1,\ldots,k}_{j_1,\ldots,j_k}}\big)\bigg] \times 
\bigg[ \dprod_{j_k=1}^{r_k} \cdots \dprod_{j_b=1}^{\rhat_b} \cdots \dprod_{j_1=1}^{r_1} Z\big(\ord{\mhat^{1,\ldots,k}_{j_1,\ldots,j_k}}\big) \bigg]$}.\]

\subsection*{The Codomain}

In the codomain of $(\cP F)^2_b$ in \cref{PFtwob}, $\Box$ is taken in $\cP X_b$.  Thus it contains the object
\begin{equation}\label{mbxbmbxbhat}
(m^b,x^b) \Box (\mhat^b, \xhat^b) = (\mtil^b, \xtil^b) \in \cP X_b
\end{equation}
with the following entries.
\[\begin{split}
\mtil^b &= (m^b, \mhat^b)\\
&= \big( m^b_1, \ldots, m^b_{r_b}, \mhat^b_1, \ldots, \mhat^b_{\rhat_b} \big) \in \cA\\
\xtil^b &= (x^b, \xhat^b)\\ 
&= \big( x^b_1, \ldots, x^b_{r_b}, \xhat^b_1, \ldots, \xhat^b_{\rhat_b} \big)\\ 
&\in (AX_b)(m^b, \mhat^b) = 
\bigg[ \prod_{j=1}^{r_b} X_b\ord{m^b_j} \bigg] \times 
\bigg[ \prod_{j=1}^{\rhat_b} X_b\ord{\mhat^b_j} \bigg].
\end{split}\]
Note that for $y \in \{m,x\}$, $\ytil^b$ has entries
\begin{equation}\label{ybarb}
\ytil^b_j = \begin{cases}
y^b_j & \text{if $1 \leq j \leq r_b$,}\\
\yhat^b_{j-r_b} & \text{if $r_b + 1 \leq j \leq r_b + \rhat_b$.}
\end{cases}
\end{equation}
Replacing the $b$-th entry in $\ang{m,x}$ with \cref{mbxbmbxbhat} yields the object
\begin{equation}\label{mxcompbmxbarb}
\begin{split}
&\ang{m,x} \circ_b (\mtil^b, \xtil^b)\\
&= \Big((m^1,x^1), \ldots, \underbrace{(\mtil^b, \xtil^b)}_{\text{$b$-th entry}}, \ldots, (m^k,x^k) \Big) \in \prod_{i=1}^k \cP X_i. 
\end{split}
\end{equation}
Similar to \cref{mhatonekjonek}, we define the following.
\begin{equation}\label{mbaronekjonek}
\begin{split}
\mtil^{1,\ldots,k}_{j_1,\ldots,j_k} 
&= m^1_{j_1} \cdots \, \mtil^b_{j_b} \cdots \, m^k_{j_k} \qquad (\text{with $1 \leq j_b \leq r_b + \rhat_b$})\\
\mtil^{1,\ldots,k} 
&= \Big\{ \cdots \, \Big\{ \cdots \, \Big\{ \mtil^{1,\ldots,k}_{j_1,\ldots,j_k} \Big\}_{j_1=1}^{r_1} \,\cdots\, \Big\}_{j_b=1}^{r_b + \rhat_b} \,\cdots\, \Big\}_{j_k=1}^{r_k} \in \cA\\
\xtil^{1,\ldots,k}_{j_1,\ldots,j_k} 
&= \big(x^1_{j_1}, \ldots, \xtil^b_{j_b}, \ldots, x^k_{j_k}\big)\\
&\in \big(X_1 \ord{m^1_{j_1}}\big) \sma \cdots \sma \big(X_b \ord{\mtil^b_{j_b}}\big) \sma \cdots \sma \big(X_k \ord{m^k_{j_k}}\big)\\
F(\xtil^{1,\ldots,k}) 
&= \Big\{ \cdots \, \Big\{ \cdots \, \Big\{ F\big(\xtil^{1,\ldots,k}_{j_1,\ldots,j_k}\big) \Big\}_{j_1=1}^{r_1} \,\cdots\, \Big\}_{j_b=1}^{r_b + \rhat_b} \,\cdots\, \Big\}_{j_k=1}^{r_k}\\
& \in \prod_{j_k=1}^{r_k} \cdots \prod_{j_b=1}^{r_b + \rhat_b} \cdots \prod_{j_1=1}^{r_1} Z\big(\ord{\mtil^{1,\ldots,k}_{j_1,\ldots,j_k}}\big) = (AZ)(\mtil^{1,\ldots,k})
\end{split}
\end{equation}
In each of the definitions of $\mtil^{1,\ldots,k}_{j_1,\ldots,j_k}$ and $\xtil^{1,\ldots,k}_{j_1,\ldots,j_k}$ above, the entry involving $\til{?}^b_{j_b}$ appears in the $b$-th spot.  By definition \cref{pfmx}, applying $\cP F$ to \cref{mxcompbmxbarb} yields the codomain of $(\cP F)^2_b$ in \cref{PFtwob}:
\begin{equation}\label{PFtwocodomain}
\begin{split}
&(\cP F)\Big( \ang{m,x} \circ_b \big( (m^b,x^b) \Box (\mhat^b, \xhat^b) \big) \Big)\\
&= \big( \mtil^{1,\ldots,k}, F(\xtil^{1,\ldots,k}) \big) \in \cP Z.
\end{split}
\end{equation}

\subsection*{The Permutation}

The first entries in the domain and the codomain of $(\cP F)^2_b$ in \cref{PFtwobdomain,PFtwocodomain} are, respectively, the objects
\begin{equation}\label{mmhatmbar}
\big( m^{1,\ldots,k}, \mhat^{1,\ldots,k} \big) \andspace \mtil^{1,\ldots,k} \in \cA.
\end{equation}
Each of these two objects is a sequence of positive integers.  By \cref{monek,mhatonekjonek,ybarb,mbaronekjonek}, there is a permutation
\begin{equation}\label{permutationsigma}
\sigma_{\{r_i\}_{i=1}^k; b, \rhat_b} \in \Sigma_{r_1 \cdots\, r_{b-1} (r_b + \rhat_b) r_{b+1} \cdots\, r_k}
\end{equation}
that takes the first object in \cref{mmhatmbar} to the second object there.  

More explicitly, if $b=k$, then
\begin{equation}\label{sigmabequalk}
\sigma_{\{r_i\}_{i=1}^k; k, \rhat_k} = \id \in \Sigma_{r_1 \cdots\, r_{k-1} (r_k + \rhat_k)}
\end{equation}
is the identity permutation.  Next suppose $1 \leq b < k$.  For $1 \leq j_i \leq r_i$ and $1 \leq i \leq k$, the permutation $\sigma_{\{r_i\}_{i=1}^k; b, \rhat_b}$ is given by
\[j_1 + \bigg[ \sum\limits_{i=2}^k r_1 \cdots\, r_{i-1} (j_i - 1) \bigg]
\mapsto j_1 + \bigg[ \sum\limits_{i=2}^k r_1 \cdots\, (r_b + \rhat_b) \cdots\, r_{i-1} (j_i - 1) \bigg].\]
On the other hand, if $1 \leq j_b \leq \rhat_b$, then  the permutation $\sigma_{\{r_i\}_{i=1}^k; b, \rhat_b}$ is given by
\[\scalebox{.95}{$
(r_1 \cdots\, r_k) + j_1 + \bigg[ \dsum\limits_{i=2}^k r_1 \cdots\, \rhat_b \cdots\, r_{i-1} (j_i - 1) \bigg] 
\mapsto j'_1 + \bigg[ \dsum\limits_{i=2}^k r_1 \cdots\, (r_b + \rhat_b) \cdots\, r_{i-1} (j'_i - 1) \bigg]$} \]
with
\[j'_i = \begin{cases}
j_i & \text{if $i \not= b$,}\\
r_b + j_b & \text{if $i = b$}.
\end{cases}\]

Thus the permutation $\sigma_{\{r_i\}_{i=1}^k; b, \rhat_b}$ in \cref{permutationsigma} defines an isomorphism
\begin{equation}\label{sigmammhatmbar}
\sigma_{\{r_i\}_{i=1}^k; b, \rhat_b} \cn \big( m^{1,\ldots,k}, \mhat^{1,\ldots,k} \big) \fto{\iso} \mtil^{1,\ldots,k} \in \cA
\end{equation}
that block permutes the unpointed finite sets
\[\begin{tikzcd}
\bigg[ \coprod\limits_{j_k=1}^{r_k} \cdots \coprod\limits_{j_1=1}^{r_1} \ufs{m^{1,\ldots,k}_{j_1,\ldots,j_k}} \bigg] \coprod \bigg[ \coprod\limits_{j_k=1}^{r_k} \cdots \coprod\limits_{j_b=1}^{\rhat_b} \cdots \coprod\limits_{j_1=1}^{r_1} \ufs{\mhat^{1,\ldots,k}_{j_1,\ldots,j_k}} \bigg] \ar{d}{\iso}\\
\bigg[ \coprod\limits_{j_k=1}^{r_k} \cdots \coprod\limits_{j_b=1}^{r_b + \rhat_b} \cdots \coprod\limits_{j_1=1}^{r_1} \ufs{\mtil^{1,\ldots,k}_{j_1,\ldots,j_k}} \bigg]
\end{tikzcd}\]
according to the permutation $\sigma_{\{r_i\}_{i=1}^k; b, \rhat_b}$.  It does not permute the elements within each unpointed finite set $\ufs{m^{1,\ldots,k}_{j_1,\ldots,j_k}}$ or $\ufs{\mhat^{1,\ldots,k}_{j_1,\ldots,j_k}}$.

The same permutation $\sigma_{\{r_i\}_{i=1}^k; b, \rhat_b}$ in \cref{permutationsigma} also takes the second entry in the domain of $(\cP F)^2_b$ in \cref{PFtwobdomain}, which is $\big( F(x^{1,\ldots,k}), F(\xhat^{1,\ldots,k})\big)$, to the second entry in the codomain of $(\cP F)^2_b$ in \cref{PFtwocodomain}, which is $F(\xtil^{1,\ldots,k})$.  Thus there is an identity morphism
\begin{equation}\label{sigmastarid}
1 \cn \big(\sigma_{\{r_i\}_{i=1}^k; b, \rhat_b}\big)_* \big( F(x^{1,\ldots,k}), F(\xhat^{1,\ldots,k})\big) \fto{=} F(\xtil^{1,\ldots,k}) \in (AZ)(\mtil^{1,\ldots,k}).
\end{equation}

Using \cref{PFtwobdomain,PFtwocodomain,sigmammhatmbar,sigmastarid}, we define the component of the $b$-th linearity constraint of $\cP F$ in \cref{PFtwob} as the isomorphism
\begin{equation}\label{PFtwobdef}
\begin{tikzcd}
\Big( \big( m^{1,\ldots,k}, \mhat^{1,\ldots,k} \big), \big( F(x^{1,\ldots,k}), F(\xhat^{1,\ldots,k})\big) \Big) \ar[shorten <=-1ex]{d}{\big( \sigma_{\{r_i\}_{i=1}^k; b, \rhat_b}, \, 1\big)}[swap]{(\cP F)^2_b \,=}\\
\big( \mtil^{1,\ldots,k}, F(\xtil^{1,\ldots,k}) \big) 
\end{tikzcd}
\end{equation}
in $\cP Z$.  The naturality of $(\cP F)^2_b$ follows from the fact that $\sigma_{\{r_i\}_{i=1}^k; b, \rhat_b}$ in \cref{sigmammhatmbar} is a block permutation, while \cref{sigmastarid} is the identity morphism.  The unity and constraint unity axioms for a $k$-linear functor follow for the same reasons and from the fact that the monoidal unit in $\cP X$ is the pair $\big( (), * \big)$ as in \cref{PXmonoidalunit}.  The constraint associativity, symmetry, and 2-by-2 axioms all follow from the fact that, between any two permuted words of the same length, there exists a unique permutation that takes one to the other.  Therefore, we have constructed a strong $k$-linear functor
\[\big( \cP F, \big\{ (\cP F)^2_i \big\}_{i=1}^k \big) \cn \prod_{i=1}^k \cP X_i \to \cP Z\]
as in \cref{prodpxpfpz}.

Note that if $k=1$, then
\[\big( \cP F, (\cP F)^2 \big) \cn \cP X_1 \to \cP Z\]
is a strict symmetric monoidal functor because, in this case, 
\[\sigma_{r_1; 1, \rhat_1} = \id \in \Sigma_{r_1 + \rhat_1}\]
is the identity permutation, as noted in \cref{sigmabequalk}.  Thus, in the case $k=1$, we recover Mandell's definition of $\cP F$ in \cite{mandell_inverseK}.

\section{Non-Symmetric Enriched Multifunctoriality}
\label{sec:Pmultifunctor}

Next we observe that the assignment
\[\begin{tikzcd}[column sep=large]
\Gacat \ar{r}{\cP} & \permcatsus
\end{tikzcd}\]
on objects \cref{mandellPX} and multimorphisms \cref{prodpxpfpz} is a non-symmetric $\Cat$-enriched multifunctor as in \cref{def:enr-multicategory-functor} with $\V = \Cat$.

\subsection*{Non-Symmetric Multifunctoriality}

The definition of $\cP$ at a multimorphism in $\Gacat$ is in \cref{pfmx,PFmorphism,PFtwobdef}.  Colored units in $\Gacat$ are identity $\Ga$-category morphisms, and $\cP$ sends them to identity symmetric monoidal functors.

\subsubsection*{Preservation of Composition: Context}

To see that $\cP$ preserves composition of multimorphisms, consider $\Ga$-category morphisms
\begin{equation}\label{FFi}
\txsma_{i=1}^k X_i \fto{F} Z \andspace \txsma_{j=1}^{n_i} W_{ij} \fto{F_i} X_i
\end{equation}
for $i \in \{1,\ldots,k\}$, with $F$ a $k$-morphism as in \cref{smashxfz} and $F_i$ an $n_i$-morphism.  We omit the comma in the subscript to simplify the notation, so $W_{ij} = W_{i,j}$.  Their composite in the multicategory $\Gacat$ is the $(n_1+\cdots+n_k)$-morphism
\begin{equation}\label{smaFiF}
\txsma_{i=1}^k \txsma_{j=1}^{n_i} W_{ij} \fto{\txsma_i F_i} 
\txsma_{i=1}^k X_i \fto{F} Z.
\end{equation}
Using \cref{Fcompnent} the composite $F \circ \big(\txsma_{i=1}^k F_i\big)$ is determined by component composite pointed functors
\begin{equation}\label{smaFiFcomponent}
\txsma_{i=1}^k \txsma_{j=1}^{n_i} W_{ij}\ord{p_{ij}} \fto{\txsma_i F_i} 
\txsma_{i=1}^k X_i\ord{p_i} \fto{F} Z\left(\ord{p_1\cdots p_k}\right)
\end{equation}
for objects $\ord{p_{ij}} \in \Fskel$, where
\[p_i = \txprod_{j=1}^{n_i}\, p_{ij} \forspace i \in \{1,\ldots,k\}.\]
The image of \cref{smaFiF} under $\cP$ is the strong $(n_1+\cdots+n_k)$-linear functor
\begin{equation}\label{PFsmaFi}
\txprod_{i=1}^k \txprod_{j=1}^{n_i} \cP W_{ij} \fto{\cP\left(F \circ (\txsma_i F_i)\right)} \cP Z.
\end{equation}

On the other hand, the images of $F$ and $F_i$ in \cref{FFi} under $\cP$ are strong multilinear functors
\[\txprod_{i=1}^k \cP X_i \fto{\cP F} \cP Z \andspace 
\txprod_{j=1}^{n_i} \cP W_{ij} \fto{\cP F_i} \cP X_i.\]
Their composite in the multicategory $\permcatsus$ is the strong $(n_1+\cdots+n_k)$-linear functor
\begin{equation}\label{prodPFiPF}
\txprod_{i=1}^k \txprod_{j=1}^{n_i} \cP W_{ij} \fto{\txprod_i\, \cP F_i} \txprod_{i=1}^k \cP X_i \fto{\cP F} \cP Z.
\end{equation}
We must show that \cref{PFsmaFi,prodPFiPF} are equal on objects, morphisms, and linearity constraints.

\subsubsection*{Preservation of Composition: Objects and Morphisms}

To check that \cref{PFsmaFi,prodPFiPF} are equal on objects, consider objects
\begin{equation}\label{mijxij}
\big(m^{ij}, x^{ij}\big) \in \cP W_{ij}
\end{equation}
for $i \in \{1,\ldots,k\}$ and $j \in \{1,\ldots,n_i\}$, with the following component objects.
\[\begin{split}
m^{ij} &= \left(m^{ij}_1, \ldots, m^{ij}_{r_{ij}}\right) \in \cA\\
x^{ij} &= \left(x^{ij}_1, \ldots, x^{ij}_{r_{ij}}\right) \in (AW_{ij})(m^{ij}) 
= \txprod_{\ell=1}^{r_{ij}} W_{ij} \ord{m^{ij}_\ell}
\end{split}\]
First we compute the image of the object
\begin{equation}\label{mijxij-tuple}
\Big\{ \Big\{ \big(m^{ij}, x^{ij}\big) \Big\}_{j=1}^{n_i} \Big\}_{i=1}^k \in \txprod_{i=1}^k \txprod_{j=1}^{n_i} \cP W_{ij}
\end{equation}
under the functor $(\cP F) \circ (\txprod_i \cP F_i)$ in \cref{prodPFiPF}.

For each $i \in \{1,\ldots,k\}$, applying $\cP F_i$ to the objects in \cref{mijxij} for $j \in \{1,\ldots,n_i\}$ yields the object
\begin{equation}\label{PFimijxij}
(\cP F_i) \left\{\big(m^{ij}, x^{ij}\big)\right\}_{j=1}^{n_i} = (m^i, x^i) \in \cP X_i
\end{equation}
with the following notation adapted from \cref{monekjonek,monek,xonekjonek}.
\[\begin{split}
m^{i1,\ldots,in_i}_{\ell_{i1}, \ldots, \ell_{in_i}} 
&= \txprod_{j=1}^{n_i} m^{ij}_{\ell_{ij}} \forspace \ell_{ij} \in \{1,\ldots,r_{ij}\}\\
m^i &= \Big\{ \,\cdots\, \Big\{ m^{i1,\ldots,in_i}_{\ell_{i1}, \ldots, \ell_{in_i}} \Big\}_{\ell_{i1}=1}^{r_{i1}} \,\cdots\, \Big\}_{\ell_{in_i}=1}^{r_{in_i}} \in \cA\\
x^{i1,\ldots,in_i}_{\ell_{i1}, \ldots, \ell_{in_i}} 
&= \left( x^{i1}_{\ell_{i1}}, \ldots, x^{in_i}_{\ell_{in_i}} \right) \in \txsma_{j=1}^{n_i} W_{ij} \ord{m^{ij}_{\ell_{ij}}}\\
x^i &= \Big\{ \,\cdots\, \Big\{ F_i\big(x^{i1,\ldots,in_i}_{\ell_{i1}, \ldots, \ell_{in_i}}\big) \Big\}_{\ell_{i1}=1}^{r_{i1}} \,\cdots\, \Big\}_{\ell_{in_i}=1}^{r_{in_i}}\\
&\in \txprod_{\ell_{in_i=1}}^{r_{in_i}} \cdots\, \txprod_{\ell_{i1}=1}^{r_{i1}} X_i \Big( \ord{m^{i1,\ldots,in_i}_{\ell_{i1}, \ldots, \ell_{in_i}}} \Big) = (AX_i)(m^i)
\end{split}\]
We omit some commas to simplify the notation, so $r_{in_i} = r_{i,n_i}$ and $\ell_{in_i} = \ell_{i,n_i}$.

Next, applying $\cP F$ to the objects in \cref{PFimijxij} for $i \in \{1,\ldots,k\}$ yields the object
\begin{equation}\label{PFmixionek}
(\cP F) \big\{ (m^i, x^i) \big\}_{i=1}^k 
= \big(m^{1,\ldots,k}, F(x^{1,\ldots,k}) \big) \in \cP Z
\end{equation}
with the following component objects.
\[\scalebox{.9}{$\begin{split}
m^{1,\ldots,k} &= \Big\{ \,\cdots\, \Big\{ \,\cdots\, \Big\{ \,\cdots\, \Big\{ \txprod_{i=1}^k m^{i1,\ldots,in_i}_{\ell_{i1}, \ldots, \ell_{in_i}} \Big\}_{\ell_{11}=1}^{r_{11}} \,\cdots\, \Big\}_{\ell_{1n_1}=1}^{r_{1n_1}} \,\cdots\, \Big\}_{\ell_{k1}=1}^{r_{k1}} \,\cdots\, \Big\}_{\ell_{kn_k}=1}^{r_{kn_k}} \in \cA\\
F(x^{1,\ldots,k}) &= \Big\{ \,\cdots\, \Big\{ \,\cdots\, \Big\{ \,\cdots\, \Big\{ 
F\big\{ F_i\big(x^{i1,\ldots,in_i}_{\ell_{i1}, \ldots, \ell_{in_i}}\big) \big\}_{i=1}^k
\Big\}_{\ell_{11}=1}^{r_{11}} \,\cdots\, \Big\}_{\ell_{1n_1}=1}^{r_{1n_1}} \,\cdots\, \Big\}_{\ell_{k1}=1}^{r_{k1}} \,\cdots\, \Big\}_{\ell_{kn_k}=1}^{r_{kn_k}}\\
&\in \txprod_{\ell_{kn_k}=1}^{r_{kn_k}} \,\cdots\, \txprod_{\ell_{k1}=1}^{r_{k1}} \,\cdots\, \txprod_{\ell_{1n_1}=1}^{r_{1n_1}} \,\cdots\, \txprod_{\ell_{11}=1}^{r_{11}} \, Z\Big( \ord{\txprod_{i=1}^k m^{i1,\ldots,in_i}_{\ell_{i1}, \ldots, \ell_{in_i}}} \Big)\\
&= (AZ)(m^{1,\ldots,k})
\end{split}$}\]
In summary, the object in \cref{PFmixionek} is the result of applying the functor $(\cP F) \circ (\txprod_i \cP F_i)$ in \cref{prodPFiPF} to the object in \cref{mijxij-tuple}.

On the other hand, for the functor in \cref{PFsmaFi}, the object equality 
\[\cP\big(F \circ (\txsma_i F_i)\big) \Big\{ \Big\{ \big(m^{ij}, x^{ij}\big) \Big\}_{j=1}^{n_i} \Big\}_{i=1}^k 
= \big(m^{1,\ldots,k}, F(x^{1,\ldots,k}) \big) \in \cP Z\]
also holds. 
\begin{itemize}
\item For the $\cA$-component, this holds because, for $\ell_{ij} \in \{1,\ldots,r_{ij}\}$, there is an equality of positive integers
\[\txprod_{i=1}^k \txprod_{j=1}^{n_i} m^{ij}_{\ell_{ij}} 
= \txprod_{i=1}^k m^{i1,\ldots,in_i}_{\ell_{i1}, \ldots, \ell_{in_i}}.\]
\item For the second component, we use the component pointed functors in \cref{smaFiFcomponent}.
\end{itemize}
This proves that \cref{PFsmaFi,prodPFiPF} are equal on objects.  A similar analysis shows that they are equal on morphisms.  Thus they are equal as functors.

\subsubsection*{Preservation of Composition: Linearity Constraints}

To see that \cref{PFsmaFi,prodPFiPF} are equal on the linearity constraints, we use the fact that, between any two permuted words of the same length, there is a unique permutation that takes one to the other, while the second factor in \cref{PFtwobdef} is an identity morphism.  Thus $\cP$ is a non-symmetric multifunctor of underlying multicategories.  The $\Cat$-enrichment of $\cP$ is discussed below.

\begin{remark}[Non-Preservation of Symmetric Group Action]\label{remark:pseudo-symmetry}
  We note that $\cP$ does not preserve the symmetric group action.  Let $(-)^\si$ denote the right action of a permutation $\si \in \Si_k$ on
  the $k$-ary operations of $\GaCat$ and $\permcatsus$, for $k > 0$.
  Then from \eqref{pfmx} we have
  \begin{align}
    (\cP F)^\si\big((m^1,x^1), \ldots , (m^k,x^k)\big)
    & = (\cP F)\big[\si\big((m^1,x^1), \ldots , (m^k,x^k)\big)\big]\nonumber\\
    & = \big(m^{\si^\inv(1),\ldots,\si^\inv(k)}, F(x^{\si^\inv(1),\ldots,\si^\inv(k)})\big)\label{eq:PFsi-1}
  \end{align}
  and
  \begin{align}
    (\cP (F^\si))\big((m^1,x^1), \ldots , (m^k,x^k)\big)
    & = \big(m^{1,\ldots,k}, F^\si(x^{1,\ldots,k})\big)\nonumber\\
    & = \big(m^{1,\ldots,k}, \pi_*^\inv F(x^{\si^\inv(1),\ldots,\si^\inv(k)})\big),\label{eq:PFsi-2}
  \end{align}
  where $\pi$ denotes the permutation of tuples that takes $m^{1,\ldots,k}$ to $m^{\si^\inv(1),\ldots,\si^\inv(k)}$ and $x^{1,\ldots,k}$ to $x^{\si^\inv(1),\ldots,\si^\inv(k)}$, induced by $\si$ permuting the rows of \cref{eq:monek-rows}.
  In \cref{eq:PFsi-2}, the second component uses the action of $\Sigma_k$ on $\Gacat$, induced by the Day convolution braiding \cref{eq:Fxi}, together with naturality of $F$ with respect to morphisms in $\Fskel$.  Then $(\pi,1)$ provides an isomorphism in $\cP Z$---generally not an identity---from 
  \cref{eq:PFsi-2} to \cref{eq:PFsi-1}.
\end{remark}

\subsection*{$\Cat$-Enrichment}

As in \cref{corollary:Ga-C-smenr} with $\C = \Cat$, the category of $k$-morphisms
\[\Gacat\mmap{Z; X_1, \ldots, X_k}\]
has
\begin{itemize}
\item as objects, $\Gamma$-category morphisms
\[F \cn \bigwedge_{i=1}^k X_i \to Z\]
as in \cref{smashxfz}, and
\item as morphisms, pointed modifications
\begin{equation}\label{ptdmodification}
\begin{tikzpicture}[xscale=2.5,yscale=2,baseline={(a.base)}]
\draw[0cell]
(0,0) node (a) {\bigwedge_{i=1}^k X_i}
(a)++(1,0) node (b) {Z}
;
\draw[1cell]  
(a) edge[bend left] node[pos=.4] {F} (b)
(a) edge[bend right] node[swap,pos=.4] {G} (b)
;
\draw[2cell] 
node[between=a and b at .47, shift={(0,0)}, rotate=-90, 2label={above,\theta}] {\Rightarrow}
;
\end{tikzpicture}
\end{equation}
between $\Gamma$-category morphisms.
\end{itemize}
Using the component pointed functors in \cref{Fcompnent}, such a pointed modification $\theta \cn F \to G$ consists of, for each $k$-tuple of objects
\begin{equation}\label{ordtup}
\ordtu{p} = \big(\ord{p_1}, \ldots, \ord{p_k}\big) \in \Fskel^{\sma k},
\end{equation} 
a pointed natural transformation
\begin{equation}\label{thetaordtup}
\begin{tikzpicture}[xscale=3.5,yscale=2,baseline={(a.base)}]
\draw[0cell]
(0,0) node (a) {\bigwedge_{i=1}^k \big(X_i \ord{p_i}\big)}
(a)++(1,0) node (b) {Z\big(\ord{p_1 \cdots\, p_k}\big)}
;
\draw[1cell]  
(a) edge[bend left] node[pos=.5] {F_{p_1 \cdots\, p_k}} (b)
(a) edge[bend right] node[swap,pos=.5] {G_{p_1 \cdots\, p_k}} (b)
;
\draw[2cell] 
node[between=a and b at .43, shift={(0,0)}, rotate=-90, 2label={above,\theta_{\ordtu{p}}}] {\Rightarrow}
;
\end{tikzpicture}
\end{equation}
that satisfies the modification axiom.  The latter encodes compatibility of $\theta$ with morphisms in $\Fskel^{\sma k}$.

On the other hand, recalling \cref{definition:permcat-multicat}, the category of $k$-morphisms
\[\permcatsus\mmap{\cP Z; \cP X_1, \ldots, \cP X_k}\]
has
\begin{itemize}
\item as objects, strong $k$-linear functors
\[\big( H, \big\{ H^2_i \big\}_{i=1}^k \big) \cn \prod_{i=1}^k \cP X_i \to \cP Z,\]
and
\item as morphisms, $k$-linear transformations
\begin{equation}\label{klineartransformation}
\begin{tikzpicture}[xscale=3.5,yscale=2,baseline={(a.base)}]
\draw[0cell]
(0,0) node (a) {\prod_{i=1}^k \cP X_i}
(a)++(1,0) node (b) {\cP Z}
;
\draw[1cell]  
(a) edge[bend left] node[pos=.4] {\big(H, \{H^2_i\}\big)} (b)
(a) edge[bend right] node[swap,pos=.4] {\big(I, \{I^2_i\}\big)} (b)
;
\draw[2cell] 
node[between=a and b at .47, shift={(0,0)}, rotate=-90, 2label={above,\vartheta}] {\Rightarrow}
;
\end{tikzpicture}
\end{equation}
between $k$-linear functors.
\end{itemize}

Given a pointed modifcation $\theta \cn F \to G$ as in \cref{ptdmodification}, its image under $\cP$ is the $k$-linear transformation
\begin{equation}\label{Pthetaklinear}
\begin{tikzpicture}[xscale=3.7,yscale=2.5,baseline={(a.base)}]
\draw[0cell]
(0,0) node (a) {\prod_{i=1}^k \cP X_i}
(a)++(1,0) node (b) {\cP Z}
;
\draw[1cell]  
(a) edge[bend left] node[pos=.45] {\big( \cP F, \{(\cP F)^2_i\}\big)} (b)
(a) edge[bend right] node[swap,pos=.45] {\big( \cP G, \{(\cP G)^2_i\}\big)} (b)
;
\draw[2cell] 
node[between=a and b at .47, shift={(0,0)}, rotate=-90, 2label={above,\cP\theta}] {\Rightarrow}
;
\end{tikzpicture}
\end{equation}
defined as follows.  If $k=0$, then the $\Gamma$-category morphism
\[F \cn J \to Z\] 
as in \cref{JFZ} is completely determined by the object $F_1(1) \in Z\ord{1}$, and similarly for $G$.  Using \cref{foneonePZ}, the 0-linear transformation $\cP \theta$ is determined by the morphism
\[\cP \theta \cn \big((1), F_1(1)\big) \to \big((1), G_1(1)\big) \in \cP Z\]
consisting of
\begin{itemize}
\item the identity morphism 
\[1 \cn (1) \to (1) \in \cA\] 
and
\item the given component morphism
\begin{equation}\label{thetaone}
\theta_1 \cn F_1(1) \to G_1(1) \in Z\ord{1}.
\end{equation}
\end{itemize} 

For $k>0$, suppose given objects 
\[(m^i,x^i) \in \cP X_i \forspace i \in \{1,\ldots,k\}\] 
as in \cref{mixi}.  With the shorthand $\ang{m,x}$ in \cref{angmx} and the definition of $(\cP F)\ang{m,x}$ in \cref{pfmx}, the component morphism
\[\begin{tikzcd}[column sep=huge]
\big( m^{1,\ldots,k}, F(x^{1,\ldots,k}) \big) \ar{r}{(\cP \theta)_{\ang{m,x}}} & 
\big( m^{1,\ldots,k}, G(x^{1,\ldots,k}) \big) \in \cP Z 
\end{tikzcd}\] 
consists of
\begin{itemize}
\item the identity morphism 
\[1 \cn m^{1,\ldots,k} \to m^{1,\ldots,k} \in \cA\]
and
\item the given component morphisms
\begin{equation}\label{thetaxonekjonek}
\begin{tikzcd}[column sep=huge]
F\big( x^{1,\ldots,k}_{j_1,\ldots,j_k} \big) \ar{r}{\theta_{x^{1,\ldots,k}_{j_1,\ldots,j_k}}} & 
G\big( x^{1,\ldots,k}_{j_1,\ldots,j_k} \big) \in Z\big( \ord{m^{1,\ldots,k}_{j_1,\ldots,j_k}} \big)
\end{tikzcd}
\end{equation}
for $1 \leq j_i \leq r_i$ and $1 \leq i \leq k$.
\end{itemize}
This finishes the definition of $\cP\theta$.

Next we check the $k$-linear transformation axioms for $\cP\theta$.
\begin{itemize}
\item The naturality of $\cP\theta$ follows from 
\begin{itemize}
\item the naturality of each component $\theta_{\ordtu{p}}$ in \cref{thetaordtup} and
\item the fact that each component of $\cP\theta$ has an identity morphism in $\cA$ as its first component.
\end{itemize}
\item If any $(m^i,x^i)$ is the monoidal unit, then $(\cP\theta)_{\ang{m,x}}$ is the identity morphism of the monoidal unit $\big((),*\big)$ in $\cP Z$ because the given $\theta$ is a \emph{pointed} modification, with each $\theta_{\ordtu{p}}$ a pointed natural transformation.
\item In the compatibility diagram \cref{eq:monoidal-in-each-variable} for $\cP\theta$ with the monoidal constraints of $\cP F$ and $\cP G$, in each of the two composites, the $\cA$-component in $\cP\theta$ is an identity morphism, while for $(\cP F)^2_b$ and $(\cP G)^2_b$ it is the same permutation in \cref{permutationsigma}.  In the other component, $(\cP F)^2_b$ and $(\cP G)^2_b$ have identity morphisms as in \cref{sigmastarid}, while $\cP\theta$ is given by the components of $\theta$ as in \cref{thetaone,thetaxonekjonek}.
\end{itemize}
Thus $\cP\theta$ is a $k$-linear transformation between $k$-linear functors.

Next we observe that the assignment
\begin{equation}\label{Pisfunctor}
\begin{tikzcd}[column sep=large]
\Gacat\mmap{Z; X_1, \ldots, X_k} \ar{r}{\cP} & \permcatsus\mmap{\cP Z; \cP X_1, \ldots, \cP X_k}
\end{tikzcd}
\end{equation}
is a functor.
\begin{itemize}
\item The assignment $\cP$ sends an identity modification $1_F$ to the identity $k$-linear transformation $1_{\cP F}$ because the component morphisms of $\theta = 1_F$ in \cref{thetaone,thetaxonekjonek} are identity morphisms.  
\item The assignment $\cP$ preserves composition of pointed modifications because, on the one hand, composition in the domain category of $\cP$ is given by vertical composition of pointed natural transformations as in \cref{thetaordtup}, which is objectwise composition in a component category of $Z$.  On the other hand, composition in the codomain category of $\cP$ is given by vertical composition of natural transformations as in \cref{klineartransformation}, which is objectwise composition in $\cP Z$.  Thus it follows from the definitions  \cref{thetaone,thetaxonekjonek} that $\cP$ preserves composition.
\end{itemize}

Non-symmetric multifunctoriality of $\cP$ with respect to multimorphisms is explained above.
The preservation of colored units and composition at the level of pointed modifications \cref{ptdmodification} both follow from the definition of $\cP\theta$ as having 
\begin{itemize}
\item an identity morphism in $\cA$ as its first component and
\item component morphisms of $\theta$, \cref{thetaone,thetaxonekjonek}, in the other component.
\end{itemize}
Thus
\[\begin{tikzcd}[column sep=large]
\Gacat \ar{r}{\cP} & \permcatsus
\end{tikzcd}\]
is a non-symmetric $\Cat$-enriched multifunctor.
This finishes the proof of \Cref{thm:invK_multifunctor}.

\begin{remark}\label{rk:redefinition}
  In the definition of the strong $k$-linear functor $\cP F$ in \cref{pfmx}, we can redefine
  \begin{enumerate}
  \item $m^{1,\ldots,k}$ in \cref{monek} to 
    \[
      m^{1,\ldots,k}_\transp = \Big\{\,\cdots\, \Big\{ m^{1,\ldots,k}_{j_1,\ldots,j_k} \Big\}_{j_k = 1}^{r_k} \,\cdots\, \Big\}_{j_1 = 1}^{r_1}
    \]
    and 
  \item $F(x^{1,\ldots,k})$ in \cref{Fxonek} to 
    \[
      F\big(x^{1,\ldots,k}_\transp\big) = \Big\{\,\cdots\, \Big\{ F\big(x^{1,\ldots,k}_{j_1,\ldots,j_k}\big) \Big\}_{j_k = 1}^{r_k} \,\cdots\, \Big\}_{j_1 = 1}^{r_1}
    \]
  \end{enumerate}
  using the lexicographic ordering instead of the reverse lexicographic ordering with respect to $j_1,\ldots,j_k$.
  Likewise, we can make a corresponding redefinition of $\cP F$ on morphisms \cref{PFmorphism} and linearity constraints \cref{PFtwobdef}.
  Call this $\cP' F$.
  Then $\cP' F$ is another strong $k$-linear functor as in \cref{prodpxpfpz}.
  Similarly, we can define $\cP' \theta$ \cref{Pthetaklinear} with the same components \cref{thetaxonekjonek} but ordered lexicographically.

  There is a canonical invertible $k$-linear transformation 
  \[
    \alpha : \cP F \fto{\iso} \cP' F 
  \]
  whose component at $\ang{ (m^i, x^i) }$ is the isomorphism in $\cP Z$
  \[\scalebox{.9}{$
    (\cP F) \ang{ (m^i, x^i) } = \left( m^{1,\ldots,k} , F(x^{1,\ldots,k}) \right) 
    \fto[\iso]{(\alpha^A, 1)}
    (\cP' F) \ang{ (m^i, x^i) } = \left( m^{1,\ldots,k}_\transp , F(x^{1,\ldots,k}_\transp )\right)$}.
  \]
  In the $\cA$-component, 
  \[
    \alpha^A : m^{1,\ldots,k} \fto{\iso} m^{1,\ldots,k}_\transp 
  \]
  is the isomorphism given by the block permutation of unpointed finite sets determined by the permutation that takes $m^{1,\ldots,k}$ to $m^{1,\ldots,k}_\transp$.  The functor
  \[
    (AZ)(\alpha^A) : (AZ)( m^{1,\ldots,k} ) \fto{\iso} (AZ)( m^{1,\ldots,k}_\transp )
  \]
  is a permutation isomorphism that takes the object $F(x^{1,\ldots,k})$ to $F(x^{1,\ldots,k}_\transp)$.  The second component in $(\alpha^A, 1)$ is the identity morphism of $F(x^{1,\ldots,k}_\transp)$.
  
  The multilinearity diagram \cref{eq:monoidal-in-each-variable} commutes in $\cP Z$ because $(\cP F)^2_b$ is componentwise of the form $(\sigma, 1)$ with $\sigma$ the block permutation of unpointed finite sets determined by the universal permutation $\sigma$ in \cref{permutationsigma}, and likewise for $(\cP' F)^2_b$.  In \cref{eq:monoidal-in-each-variable} the first components of the two composites are equal as morphisms in $\cA$ by uniqueness of such permutations.  The second component is an identity morphism for each composite.  With these consistent changes, $\cP'$  is another extension of the functor $\cP$ to a non-symmetric Cat-enriched multifunctor.
\end{remark}

\section{Application: Inverse \texorpdfstring{$K$}{K}-Theory Preserves Associative Algebras}
\label{sec:preservation}

As a consequence of \Cref{thm:invK_multifunctor}, the $\Cat$-enriched inverse $K$-theory $\cP$ preserves algebraic structures parametrized by non-symmetric $\Cat$-enriched multicategories, in particular, the $\Cat$-enriched associative operad.

We begin with a description of $E_1$-algebras in permutative categories.
\begin{definition}[\cite{elmendorf-mandell}]\label{def:ringcat}
A \emph{ring category} is a tuple
\[\big(\C,(\oplus,\zero,\xiplus),(\otimes,\tu),(\fal,\far)\big)\]
consisting of the following data.
\begin{description}
\item[The Additive Structure] 
$(\C,\oplus,\zero,\xiplus)$ is a permutative category. 
\item[The Multiplicative Structure] 
$(\C,\otimes,\tu)$ is a strict monoidal category.
\item[The Factorization Morphisms] 
$\fal$ and $\far$ are natural transformations
\begin{equation}\label{ringcatfactorization}
\begin{tikzcd}[row sep=tiny,column sep=huge]
(A \otimes C) \oplus (B \otimes C) \ar{r}{\fal_{A,B,C}} & (A \oplus B) \otimes C\\
(A \otimes B) \oplus (A \otimes C) \ar{r}{\far_{A,B,C}} & A \otimes (B \oplus C)
\end{tikzcd}
\end{equation}
for objects $A,B,C\in\C$, which are called the \emph{left factorization morphism} and the \emph{right factorization morphism}, respectively.
\end{description}
We often abbreviate $\otimes$ to concatenation, with $\tensor$ always taking precedence over $\oplus$ in the absence of clarifying parentheses.  The subscripts in $\xiplus$, $\fal$, and $\far$ are sometimes omitted. 

The above data are required to satisfy the following seven axioms for all objects $A$, $A'$, $A''$, $B$, $B'$, $B''$, $C$, and $C'$ in $\C$.  Each diagram is required to be commutative.
\begin{description}
\item[The Multiplicative Zero Axiom] 
\[\begin{tikzcd}[column sep=large]
\boldone \times \C \ar{d}[swap]{\zero \times 1_{\C}} \ar{r}{\iso} & \C \ar{d}{\zero} & \C \times \boldone \ar{l}[swap]{\iso} \ar{d}{1_{\C} \times \zero}\\
\C \times \C \ar{r}{\otimes} & \C & \C \times \C \ar{l}[swap]{\otimes}
\end{tikzcd}\]
In this diagram, the top horizontal isomorphisms drop the $\boldone$ argument.  Each $\zero$ denotes the constant functor at $\zero \in \C$ and $1_\zero$.
\item[The Zero Factorization Axiom] 
\[\begin{aligned}
\fal_{\zero,B,C} &= 1_{B \otimes C} &\qquad \far_{\zero,B,C} &= 1_{\zero}\\
\fal_{A,\zero,C} &= 1_{A \otimes C} & \far_{A,\zero,C} &= 1_{A \otimes C}\\
\fal_{A,B,\zero} &= 1_{\zero} & \far_{A,B,\zero} &= 1_{A \otimes B}
\end{aligned}\]
\item[The Unit Factorization Axiom] 
\[\begin{split}
\fal_{A,B,\tu} &= 1_{A \oplus B}\\
\far_{\tu,B,C} &= 1_{B \oplus C}
\end{split}\]
\item[The Symmetry Factorization Axiom] 
\[\begin{tikzcd}
AC \oplus BC \ar{d}[swap]{\xiplus} \ar{r}{\fal} & (A \oplus B)C \ar{d}{\xiplus 1_C}\\
BC \oplus AC \ar{r}{\fal} & (B \oplus A)C
\end{tikzcd}\qquad
\begin{tikzcd}
AB \oplus AC \ar{d}[swap]{\xiplus} \ar{r}{\far} & A(B \oplus C) \ar{d}{1_A \xiplus}\\
AC \oplus AB \ar{r}{\far} & A(C \oplus B)
\end{tikzcd}\]
\item[The Internal Factorization Axiom]  
\[\begin{tikzcd}[cells={nodes={scale=.8}}]
AB \oplus A'B \oplus A''B \ar{d}[swap]{1 \oplus \fal} \ar{r}{\fal \oplus 1} & (A \oplus A')B \oplus A''B \ar{d}{\fal}\\
AB \oplus (A' \oplus A'')B \ar{r}{\fal} & (A \oplus A' \oplus A'')B
\end{tikzcd}\qquad
\begin{tikzcd}[cells={nodes={scale=.8}}]
AB \oplus AB' \oplus AB'' \ar{d}[swap]{1 \oplus \far} \ar{r}{\far \oplus 1} & A(B \oplus B') \oplus AB'' \ar{d}{\far}\\
AB \oplus A(B' \oplus B'') \ar{r}{\far} & A(B \oplus B' \oplus B")
\end{tikzcd}\]
\item[The External Factorization Axiom] 
\[\begin{tikzcd}[column sep=large,cells={nodes={scale=.8}}]
ABC \oplus A'BC \ar{d}[swap]{\fal_{AB,A'B,C}} \ar{r}{\fal_{A,A',BC}} & (A \oplus A')BC \ar[equal]{d}\\
(AB \oplus A'B)C \ar{r}{\fal_{A,A',B} 1_C} & (A \oplus A')BC
\end{tikzcd}
\qquad
\begin{tikzcd}[column sep=large,cells={nodes={scale=.8}}]
ABC \oplus AB'C \ar{d}[swap]{\far_{A,BC,B'C}} \ar{r}{\fal_{AB,AB',C}} & (AB \oplus AB')C \ar{d}{\far 1_C}\\
A(BC \oplus B'C) \ar{r}{1_A \fal_{B,B',C}} & A(B \oplus B')C
\end{tikzcd}\]
\[\begin{tikzcd}[column sep=huge]
ABC \oplus ABC' \ar{d}[swap]{\far_{A,BC,BC'}} \ar{r}{\far_{AB,C,C'}} & AB(C \oplus C') \ar[equal]{d}\\
A(BC \oplus BC') \ar{r}{1_A \far_{B,C,C'}} & AB(C \oplus C')
\end{tikzcd}\]
\item[The 2-By-2 Factorization Axiom]
\[\begin{tikzpicture}[xscale=3,yscale=1,baseline={(x2.base)}]
\def\h{1} \def\v{-1}
\draw[0cell=.8] 
(0,0) node (x0) {A(B \oplus B') \oplus A'(B \oplus B')}
(x0)++(-\h,\v) node (x1) {AB \oplus AB' \oplus A'B \oplus A'B'}
(x1)++(1.3*\h,\v) node (x2) {(A \oplus A')(B \oplus B')}
(x1)++(0,2*\v) node (x3) {AB \oplus A'B \oplus AB' \oplus A'B'}
(x3)++(\h,\v) node (x4) {(A \oplus A')B \oplus (A \oplus A')B'}
;
\draw[1cell=.9] 
(x1) edge node[pos=.4] {\far \oplus \far} (x0)
(x0) edge node {\fal} (x2)
(x1) edge node[swap] {1 \oplus \xiplus \oplus 1} (x3)
(x3) edge node[swap,pos=.4] {\fal \oplus \fal} (x4)
(x4) edge node[swap] {\far} (x2)
; 
\end{tikzpicture}\]
\end{description}
This finishes the definition of a ring category.  Moreover, a ring category is said to be \emph{tight} if $\fal$ and $\far$ in \cref{ringcatfactorization} are natural isomorphisms.
\end{definition}

Let $\As$ denote the $\Cat$-enriched associative operad, as in \cite[11.1.1]{cerberusIII}.  By \cite[11.2.16]{cerberusIII}, $\Cat$-enriched multifunctors from $\As$ to $\permcatsu$ parametrize ring categories.
Since $\As$ is the free symmetric operad on the terminal non-symmetric operad, which we denote $\As'$, ring categories are also parametrized by non-symmetric $\Cat$-enriched multifunctors
\[
\As' \to \permcatsu.
\]
Since the multilinear functors in $\permcatsus$ are \emph{strong}, $\As'$-algebras in $\permcatsus$ have invertible factorization morphisms
and are, therefore, tight ring categories.
Combining this characterization with \Cref{thm:invK_multifunctor} yields the following result.

\begin{corollary}\label{Kinven}
The non-symmetric $\Cat$-enriched inverse $K$-theory multifunctor
\[\begin{tikzcd}[column sep=large]
\Gacat \ar{r}{\cP} & \permcatsus
\end{tikzcd}\]
in \Cref{thm:invK_multifunctor} sends each monoid in $\Gacat$ to a tight ring category (\Cref{def:ringcat}).
\end{corollary}
\noindent Thus, tight ring categories arise naturally as the outputs of the $\Cat$-enriched inverse $K$-theory, $\cP$.

\bibliographystyle{sty/amsalpha3}
\bibliography{references}

\providecommand{\bysame}{\leavevmode\hbox to3em{\hrulefill}\thinspace}
\providecommand{\MR}{\relax\ifhmode\unskip\space\fi MR }
\providecommand{\MRhref}[2]{%
  \href{http://www.ams.org/mathscinet-getitem?mr=#1}{#2}
}
\providecommand{\nopubyear}{$\infty$}
\providecommand{\doi}[1]{%
  doi:\href{https://dx.doi.org/#1}{\nolinkurl{#1}}}
\providecommand{\arxiv}[1]{%
  arXiv:\href{https://arxiv.org/abs/#1}{#1}}
\begin{thebibliography}{Elm$\infty$AA}

\bibitem[Cis99]{cisinski-thomason-correction}
D.-C. Cisinski, \emph{La classe des morphismes de {D}wyer n'est pas stable par
  retractes}, Cahiers Topologie G\'{e}om. Diff\'{e}rentielle Cat\'{e}g.
  \textbf{40} (1999), no.~3, 227--231.

\bibitem[Elm$\infty$]{elmendorf-multi-inverse}
A.~D. Elmendorf, \emph{Multiplicativity in {M}andell’s inverse {$K$}-theory},
  Submitted. \arxiv{2110.07512}

\bibitem[EM06]{elmendorf-mandell}
A.~D. Elmendorf and M.~A. Mandell, \emph{Rings, modules, and algebras in
  infinite loop space theory}, Adv. Math. \textbf{205} (2006), no.~1, 163--228.
  \doi{10.1016/j.aim.2005.07.007}

\bibitem[EM09]{elmendorf-mandell-perm}
\bysame, \emph{Permutative categories, multicategories and algebraic
  {$K$}-theory}, Algebr. Geom. Topol. \textbf{9} (2009), no.~4, 2391--2441.
  \doi{10.2140/agt.2009.9.2391}

\bibitem[GJO17]{gjo1}
N.~Gurski, N.~Johnson, and A.~M. Osorno, \emph{{$K$}-theory for 2-categories},
  Adv. Math. \textbf{322} (2017), 378--472. \doi{10.1016/j.aim.2017.10.011}

\bibitem[JY$\infty$]{cerberusIII}
N.~Johnson and D.~Yau, \emph{{B}imonoidal {C}ategories, {$E_n$}-{M}onoidal
  {C}ategories, and {A}lgebraic {$K$}-{T}heory. {V}olume {III}: {F}rom
  {C}ategories to {S}tructured {R}ing {S}pectra}, available at
  \url{https://nilesjohnson.net}.

\bibitem[JY21]{johnson-yau}
\bysame, \emph{{2}-{D}imensional {C}ategories}, Oxford University Press, New
  York, 2021. \doi{10.1093/oso/9780198871378.001.0001}

\bibitem[{Man}10]{mandell_inverseK}
M.~A. {Mandell}, \emph{{An inverse {$K$}-theory functor}}, {Doc. Math.}
  \textbf{15} (2010), 765--791.

\bibitem[Seg74]{segal}
G.~Segal, \emph{Categories and cohomology theories}, Topology \textbf{13}
  (1974), 293--312. \doi{10.1016/0040-9383(74)90022-6}

\bibitem[Tho80]{thomason-model-str}
R.~W. Thomason, \emph{Cat as a closed model category}, Cahiers Topologie
  G\'{e}om. Diff\'{e}rentielle \textbf{21} (1980), no.~3, 305--324.

\bibitem[Tho95]{thomason}
\bysame, \emph{Symmetric monoidal categories model all connective spectra},
  Theory Appl. Categ. \textbf{1} (1995), No. 5, 78--118.

\bibitem[Yau20]{yau-hqft}
D.~Yau, \emph{Homotopical quantum field theory}, World Scientific, Singapore,
  2020.

\bibitem[YJ15]{bluemonster}
D.~Yau and M.~W. Johnson, \emph{A foundation for {PROP}s, algebras, and
  modules}, Mathematical Surveys and Monographs, vol. 203, American
  Mathematical Society, Providence, RI, 2015. \doi{10.1090/surv/203}

\end{thebibliography}
\end{document}